\documentclass{compositio}
\usepackage{amsmath,amsfonts,amssymb,mathptmx}
\usepackage{thesis_pub,amscd,amsthm}
\usepackage[all]{xy}

\newtheorem{theorem}{Theorem}[section] 
\newtheorem{lemma}[theorem]{Lemma}   \newtheorem{proposition}[theorem]{Proposition} 
\newtheorem{corollary}[theorem]{Corollary}
{\theoremstyle{definition}  \newtheorem{definition}[theorem]{}  }
{\theoremstyle{remark} \newtheorem{remark}[theorem]{Remark}  }

\begin{document}

\title{On the arithmetic of twists of superelliptic curves}
\author{Sungkon Chang}
\email{changsun@mail.armstrong.edu}
\address{Department of Mathematics\\ Armstrong Atlantic State University \\
			11935 Abercorn St \\ Savannah, GA 31419}

\begin{abstract} In this paper, we consider a family of twists of a superelliptic curve over a global field, and obtain results on the distribution of the Mordell-Weil rank of these twists. Our results have applications to the distribution  of the number of rational points. 
\end{abstract}

\maketitle

\section{Introduction} \label{sec:intro}

By Faltings' theorem, a (smooth complete geometrically irreducible)  curve over a number field has finitely many rational points.
By \cite{chm:1997}, it is widely believed that the number of rational points of a curve of genus $>1$ over a number field $K$ is bounded in terms of the genus of the curve. In \cite{mazur:1986},  prior to \cite{chm:1997}, Mazur asked whether the number of rational points can be bounded in terms of the genus and the Mordell-Weil rank of its Jacobian variety. For the case of twists of curves, in \cite{silverman:1993}, Silverman proves that Mazur's question has  a positive answer.  However, for general cases, this question is totally open. 

By Silverman's result, given a curve of genus $>1$ over a number field, finding infinitely many twists with a bounded number of rational points becomes a problem of finding infinitely many twists with bounded Mordell-Weil rank.  Even for special cases such as Thue equations (see \cite{dino:2002}), an answer to this problem is sometimes not known.  For the case of elliptic curves, by  Kolyvagin's result \cite{kolyvagin:1988} and the modularity of elliptic curves proved by Wiles et al, 
results such as \cite{ono:1998}, in which quadratic twists with analytic rank $0$ are computed, imply that given an elliptic curve over $\ratn$, there are infinitely many quadratic twists with Mordell-Weil rank $0$, i.e., algebraic rank $0$. 
There are also  results of this type such as Heath-Brown's \cite{heath-brown:1994}, \cite{wong:1999} and \cite{yu:2003} which rather directly show that there is a \lq\lq positive proportion\rq\rq\ of algebraic rank-$0$ quadratic twists of certain elliptic curves.

In this paper, we consider a family of twists of  superelliptic curves over a global field, and prove that for these twists, the problem of finding infinitely many twists with bounded Mordell-Weil rank has a positive answer and, hence, there are infinitely many twists with bounded number of rational points if the genus $>1$. Some Thue equations can be mapped down to superelliptic curves considered in this paper and, hence, for these Thue equations, this problem has a positive answer. For the case of superelliptic curves over a constant field and the case of hyperelliptic curves over $\ratn$, finer results are obtained.  Especially,  using superelliptic curves over a constant field, we show that there are (infinitely many twists of) curves of arbitrarily large genus over a function field with Mordell-Weil rank $0$.  Over a number field, examples of such curves are not known.
In this section, we also introduce the application of our result  to cubic twists of some elliptic curves.

Let $n \ge 2$ be a positive integer, and let $R$ be an integral domain of characteristic not dividing $n$, with field of fractions $K$.
Let $f(x)$ be a monic polynomial in $ R[x]$ such that $n$ is coprime to $\deg(f)$, and $f(x)$ has distinct roots.
In this paper, \textit{a superelliptic curve} is the  projective $K$-model of the affine plane curve $ y^n = f(x)$.

\begin{definition}\label{def:CDJD}
Let $K$ be a field, and let $\ell$ be a prime number different from $\Char K$. Let $C/K$ be the normalization
of a \super\ given by  $y^\ell=f(x)$. For $D \in K^*$, we denote by $C_D/K$ the normalization of the curve given by $y^\ell = D^d\, f(x/D)$ where $d:=\deg(f)$, and by $J_D/K$, the Jacobian variety of $C_D$.
         The Jacobian variety  $J_D$ is called \textit{an $\ell$-th power twist of $J$}.
For the case of hyperelliptic curves (where $\ell=2$), the plane curve $D\, y^2 = f(x)$ is isomorphic to  $y^2 = D^d\, f(x/D)$.        
We denote by $rank\ J_D(K)$ the Mordell-Weil rank of $J_D(K)$ if $J_D(K)$ is a finitely generated (abelian) group.
Let $\zetal$ be a primitive $\ell$-th root of unity, let $F$ be $K(\zetal)$, and let $\lambda:=1- \zetal$.  
Throughout the paper, we denote also by $\lambda$ the endomorphism $1 - \zetal$ on $J_F$ defined in \cite{schaefer:1998}, Sec 3, and by $\sellam(J,F)$ \textit{the $\lambda$-Selmer group of $J_F$}.
\end{definition}
\par
In this work, we shall consider both the number field case and the function field case.  We begin by introducing our results for number fields.
Let $k>1$ be a positive integer, and let us denote by $\curlyP_k(X)$ the set of all positive $k$-th power-free integers up to $X$.  
Given a polynomial $f(x)$,  let $\Delta_f$ denote the discriminant of $f(x)$.
 Let $F$ be the field of fractions of a Dedekind domain $\OF$, and let $\curlyD$ be a set of prime ideals of $\OF$.  A nonzero  element $D$ of $\OF$ is \textit{supported by $\curlyD$} if $D\OF$ is divisible only by prime ideals contained in $\curlyD$. 

  \medskip
\copytheorem{Theorem \ref{MainTheorem1}}
Let $K$ denote the $\ell$-th cyclotomic field extension $\ratn(\zetal)$ where $\ell$ is a regular prime number.
Let $f(x)$ be a monic polynomial of prime degree $p$ defined over $\zz$ such that $f(x)$ is  irreducible over $K$, and $\ell \ne p$.
Let $C/\ratn$  be the normalization of the \super\ $y^\ell=f(x)$, and let $J/\ratn$ be the Jacobian variety of $C$.
\par
Let $D_0$ be a positive integer.  Let $N:=\dimp\ell\sellam(J_{D_0},K)$, and let $M$ be the number of prime ideals of $\OK$ dividing $\ell \Delta_f\, D_0$.
Then there is a set $\curlyD$ of prime numbers with  Dirichlet density at least $(p-1)/\big( \ell^{(N+M+1)p!} (\ell-1) p!  \big)$ such that whenever a positive integer $D$ is supported by $\curlyD$,  
\begin{equation*}
 \dimp\ell\sellam(J_{D_0\, D}, K) = N.
               \end{equation*}
Moreover,  there is a positive constant $\ep<1$ depending on $C$ and $D_0$ such that 
\begin{equation*}
\Sharp\ \set{ \DPlX :\dimp\ell \sellam(J_D,K)=N} \gg_{J,D_0} \frac{ X }{ (\log X)^{\ep} }.
\end{equation*}
\rm
\par\noindent
This theorem is proved using Schaefer's description of the $\lambda$-Selmer group \cite{schaefer:1998}.
The case of superelliptic curves $y^\ell = x^2 -A$ is considered by Stoll in \cite{stoll:1998}, which in fact inspired our work (see Corollary \ref{cor:cubic} on p.\pageref{cor:cubic}).
Schaefer's description is more complicated when $\ell \mid \deg(f)$; see \cite{poonen:1997}.
\par
Finer results on the distribution of Selmer ranks of twists of the Jacobian variety of a curve have only been obtained in very special cases.
Let $T(X)$ be the set of nonzero square-free integers $D$ such that $\abs{D} < X$. 
In \cite{heath-brown:1994}, Heath-Brown studies  the distribution of $2$-Selmer ranks of quadratic twists of the elliptic curve $E/\ratn$ given by $y^2=x^3-x$.  For example, given a positive integer $n$,  he computes explicitly
$$ \lim_{X\to\infty}\ 
        \frac{ \Sharp\,\set{D \in T(X) : \dimp2\Sel2(E_D,\ratn)=n,\ D\equiv 1 \mod 8} }
                { \Sharp\,\set{D \in T(X) :  D\equiv 1 \mod 8} }.$$
        In \cite{yu:2003}, for infinitely many elliptic curves $E/\ratn$ with a nontrivial rational $2$-torsion point, Gang Yu computes an upper bound for a certain average of $\dimp2\Sel2(E_D,\ratn)$.
        Let $K:=\ratn(\zeta_3)$, and let $E/\ratn$ be the elliptic curve given by $y^2=x^3-A$ where $A$ is an integer.  Let $\lambda$ be the endomorphism defined in \ref{def:CDJD} for $\ell=3$.       
In \cite{chang:2004}, we computed  an upper bound for a certain average of $\lambda$-Selmer rank of $E_D$ over $K$ where $E_D$ is a quadratic twist of $E/\ratn$.
Given an integer $n$,  computed also in \cite{chang:2004} is a lower bound for 
$$\underset{X\to\infty}{ \lim\inf}\ %
        \frac{ \Sharp\,\set{D \in T(X) : \dimp3\sellam(E_D,K)\le 2n,\ D\equiv 1 \mod 12A} }
                { \Sharp\,\set{D \in T(X) :  D\equiv 1 \mod 12A} }.$$
The results in \cite{heath-brown:1994}, \cite{yu:2003}, and \cite{chang:2004} have applications to the distribution of Mordell-Weil ranks of quadratic twists of elliptic curves considered in these papers.            
Using Theorem \ref{MainTheorem1}, we obtain the following result on the distribution of Mordell-Weil ranks of  $\ell$-th power twists of the Jacobian variety of $C$.

\medskip
\copytheorem{Corollary \ref{cor:rank<=n}}
Assume the same hypotheses in Theorem \ref{MainTheorem1}.  
Then there is a positive constant $\ep<1$ such that 
\begin{equation}
\Sharp\ \set{ \DPlX :\rank J_D(\ratn)\le N} \gg_{J,D_0} \frac{ X }{ (\log X)^{\ep} }.\label{eq:MT3}
\end{equation}
\rm
If a Jacobian variety $J$ considered in Corollary \ref{cor:rank<=n} has an $\ell$-th power twist with $\lambda$-Selmer rank $0$, then the corollary implies that there are infinitely many $\ell$-th power twists with Mordell-Weil rank $0$. 

\medskip
\copytheorem{Theorem \ref{MainTheorem3}}
Let $K$ denote the $\ell$-th cyclotomic field extension $\ratn(\zetal)$ where $\ell$ is a regular prime number.
Let $f(x)$ be a monic polynomial   defined over $\zz$ such that $f(x)$ has a root in $K$, and $\ell \nmid \deg(f)$.
Let $C/\ratn$  be the normalization of the \super\ $y^\ell=f(x)$.
\par
Given a positive integer $n$, there is a positive constant $\ep<1$ depending on $C$ and $n$ such that 
\begin{equation*}
\Sharp\,\set{ \DPlX : \dimp\ell\Sel{\lambda}(J_D, K) > n }
                \gg_{C,n}\frac{X}{(\log X)^\ep}.
\end{equation*}
In particular,
$
\lim\sup_D\ \dimp\ell\Sel{\lambda}(J_D, K) =\infty.
$

\medskip\rm
\indent
Recall that for $D \in \ratn^*$,
$\Sel2(J_D,\ratn) \cong J_D(\ratn) / 2 J_D(\ratn) \oplus \sha(J_D,\ratn)[2]$
where $\sha(J_D,\ratn)$ is the Tate-Shafarevich group of $J_D/\ratn$.
Suppose that $\ell=2$, and that $f(x)$ is any polynomial of odd degree with a root in $\ratn$ such that $f(x)$ has distinct roots. Theorem \ref{MainTheorem3} implies in particular that the $2$-Selmer groups of quadratic twists of the hyperelliptic curve $y^2=f(x)$ can be arbitrarily large, and this result seems to be new. However, for some elliptic curves, more is known. For instance, it is proved in \cite{atake:2001} as a generalization of Lemmermeyer's work \cite{lemm:1998} that the $2$-part of the Tate-Shafarevich groups of quadratic twists of the elliptic curves considered in the paper can be arbitrarily large. 
\par
Consider a Fermat curve $x^\ell + y^\ell = 1$.  It is in fact isomorphic to a superelliptic curve $y^\ell = f(x)$ (where $\ell \nmid \deg(f)$), and using Theorem \ref{MainTheorem3}, we obtain the following result on Fermat twists.

\medskip
\copytheorem{Corollary \ref{cor:fermat}}
Let $K$ denote the $\ell$-th cyclotomic field extension $\ratn(\zetal)$ where $\ell$ is an odd regular prime number.
Let $F_D/K$ be the Fermat curve given by $x^\ell + y^\ell =D\, z^\ell$ where $D\in \zz$ is nonzero, and let $\Jac(F_D)$ be the Jacobian variety of $F_D/K$. Let $\zetal$ denote the automorphism of order $\ell$ on $F_D$ given by $x\mapsto x$, $y \mapsto y$, and $z \mapsto z\, \zetal$.  Let $\lambda$ be the endomorphism $1-[\zetal]$ on $\Jac(F_D)$ where $[\zetal]$ denotes the automorphism on $\Jac(F_D)$ induced by the automorphism on $F_D$.  Then
$$ \lim\sup_D\ \dimp\ell\sellam(\Jac(F_D),K) =\infty.$$
\rm
\indent
Let us further consider applications of Theorem \ref{MainTheorem1}.
Let $C/\ratn$ be the normalization of a superelliptic curve of genus $>1$, and $C_D$, the twist of $C$ as defined in \ref{def:CDJD} for some $D \in \ratn^*$. 
We apply Silverman's result \cite{silverman:1993}, p.~234, that if $K$ is a number field, then the number of $K$-rational points on a twist $C_D$ of $C/K$ is  bounded in terms of a constant $\gamma$ depending on $C$, and in terms of the Mordell-Weil rank of the Jacobian variety $J_D$, namely, 
\begin{equation}\label{eq:sil-bound}
 \Sharp\, C_D(K) < \gamma\  7^{\rank J_D(K)}.
 \end{equation}
\begin{corollary}\label{cor:sols}
Let $C$ be the normalization of a \super\ considered in Theorem \ref{MainTheorem1}. Let $N:= \dimp\ell\sellam(J_{D_0},K)$ for a nonzero positive integer $D_0$.  If the genus of $C$ is $>1$, then there are positive constants $\ep <1$, and $\gamma$ depending on $C$ and $D_0$ such that 
\begin{equation*} 
\Sharp\ \set{ \DPlX : \Sharp\ C_D(\ratn) \le \gamma\ 7^{N}} \gg  \frac{ X }{ (\log X)^{\ep} }.
\end{equation*}
\end{corollary}
\begin{proof}
The proof is left to the reader.\end{proof}
For the case of hyperelliptic curves $C$, using \cite{stoll:2004}, Theorem 1, we obtain a sharper upper bound on the number of rational points on quadratic twists of $C/K$. 

\medskip
\copytheorem{Corollary \ref{cor:stoll}}
Let $C/\ratn$ be the normalization of a hyperelliptic curve $y^2=f(x)$ where $f(x)\in\zz[x]$ is monic and irreducible over $\ratn$, and has odd prime degree $p\ge 5$.
Suppose that there is a positive integer $D_0$ such that $N:=\dimp2\Sel{2}(J_{D_0},\ratn) < (p-1)/2$. 
Then there is a positive constant $\ep<1$ depending on $C$ and $D_0$  such that 
\begin{equation}
\Sharp\ \set{ \DPtwoX : \Sharp\ C_D(\ratn)\le 2N+1 }
                \gg \frac{X}{(\log X)^\ep}.
\end{equation}
\rm
\indent
J.~Silverman asked (see  \cite{ono:1998}, p.~653.) whether given an elliptic curve $E/\ratn$,  there are infinitely many prime numbers $p$ for which either $E_p$ or $E_{-p}$ has Mordell-Weil rank zero.
We can show:
\begin{corollary}\label{cor:sil-ques}
Let $E/\ratn$ be an elliptic curve without $\ratn$-rational $2$-torsion  points.
If $\dimp2\linebreak[0]\Sel{2}(E,\ratn) =0$, then 
there is a set $\curlyD$ of prime numbers  with positive Dirichlet density such that $\rank E_p(\ratn)=0$ for all $p\in \curlyD$.
In particular, there are infinitely many prime numbers $p$ such that $\rank E_p(\ratn)=0$.
\end{corollary}
\begin{proof}
The proof is left to the reader.\end{proof}
In \cite{ono:1998}, Corollary 3,  Ono and Skinner proved that  the question of Silverman has a positive answer for all elliptic curves with conductor $\le 100$. Theorem \ref{MainTheorem1} and Corollary \ref{cor:sil-ques} in particular imply that there are infinitely many elliptic curves over $\ratn$ for which the question of Silverman has a positive answer: Let $E/\ratn$ be the elliptic curve $y^2=x^3-2$. Then,  $\SEL(E,\ratn)=0$ and, by Theorem \ref{MainTheorem1}, there are infinitely many quadratic twists $E_D$ with $\SEL(E_D,\ratn)=0$.  Now, Corollary \ref{cor:sil-ques} implies a positive answer to the question of Silverman for these elliptic curves $E_D$.
\par
Little is known about the distribution of quadratic twists of an elliptic curve with Mordell-Weil rank $1$.  Vatsal's result \cite{vatsal:1998} is unconditional: He proved   that for the elliptic curve $E= X_0(19)$, 
\begin{equation*}
\Sharp\, \set{ \abs{D} < X : \rank E_D(\ratn) = 1 } \gg X.
\end{equation*}
Assuming the Riemann Hypothesis, Iwaniek and Sarnak proved in \cite{iwaniec-sarnak:2000} that 
$\Sharp\, \set{\abs{D} < X : \linebreak \rank E_D(\ratn) = 1 } \gg_{E} X$  for all elliptic curves $E/\ratn$.
We prove here:
\begin{corollary}
Assume the finiteness of  the Tate-Shafarevich groups  of all elliptic curves over $\ratn$. 
Let $E/\ratn$ be an elliptic curve without $\ratn$-rational $2$-torsion points such that 
$\dimp2\Sel{2}(E_{D_0},\ratn)=1$ for some positive integer $D_0$.
Then there is a positive constant $\ep <1$ such that 
\begin{equation}\label{eq:bound:rank1}
\Sharp\ \set{ D \in \curlyP_2(X) : \rank E_D(\ratn)=1,\ \sha(E_D,\ratn)[2]=\set{0}} \gg_{E,D_0} \frac{X}{(\log X)^{\ep}}.
\end{equation}
\end{corollary}
\begin{proof} 
Let $E$ be an elliptic curve with $\SEL(E,\ratn)=1$.
 The finiteness of the Tate-Shafarevich group $\sha(E,\ratn)$ of $E/\ratn$ implies that the Cassels-Tate pairing 
$ \sha(E,\ratn) \times \sha(E,\ratn) \to \ratn / \zz$
is a non-degenerate alternating bilinear pairing.  By this pairing,  $\dimp2\sha(E,\ratn)[2] \le 1$ implies 
$\rank E(\ratn)=1$ and  $\dimp2\sha(E,\ratn)[2]=0$.
Then (\ref{eq:bound:rank1}) follows immediately from Theorem \ref{MainTheorem1}.\end{proof}
 \par
Very little is known about the distribution of Mordell-Weil ranks of cubic twists of an elliptic curve.
Let $E/\ratn$ be the elliptic curve given by $x^3+y^3=1$. 
 In \cite{lieman:1994}, Lieman showed that given an integer $c$ and a prime number $p$, there are infinitely many cubic twists  $E_D : x^3+y^3=D$  such that $D\equiv c \mod p$ and $\rank E_D(\ratn)=0$.  
 We show here
 
 \medskip
\copytheorem{Corollary \ref{cor:cubic}}
Let $E/\ratn$ be an elliptic curve given by $y^2=x^3-A$ where $A$ is a positive square-free integer such that 
$A\equiv 1$ or $25 \mod 36$ and $\dimp3 \Cl(\ratn(\sqrt{-A}))[3]=0$.  For a non-zero cube-free integer $D$, let $E_D$ be the cubic twist: $y^2=x^3-A\,D^2$.  Then there is a positive integer $\ep <1$ such that 
\begin{equation}\label{eq:cubic}
 \Sharp\, \set{ D \in \curlyP_3(X) : \rank E_D(\ratn) = 0 } \gg \frac{ X }{ (\log X)^\ep }.
 \end{equation}
\rm
\indent
  \textit{A function field of one variable over an arbitrary field $k$} is a field extension $K$ of $k$ with transcendence degree $1$ such that $K$ is finitely generated over $k$, and $k$ is algebraically closed in $K$. 
\textit{A \globalfieldv} is  a function field $K$ of one variable over a finite field $k$ with a non-archimedean absolute value $v_\infty$ on $K/k$ of \textit{degree $1$}. 
Such function fields correspond to smooth complete curves $\curlyZ$ over $k$ with a $k$-rational point $p_\infty$  corresponding to the absolute value $v_\infty$.  For this type of function fields, we choose $\OK:=\set{ \al \in K : \abs{\al}_{\vinf} > 1}$, as a ring of integers in $K$, i.e., $\OK$ is the ring of regular functions on the open subset $\curlyZ\minuS\set{p_\infty}$.
In Corollary \ref{cor:hyper-FF-no2} and Theorem \ref{MainTheorem3'}, we prove function field analogues of Theorem \ref{MainTheorem1} and \ref{MainTheorem3}.

\begin{definition}\label{def:ff-setup}
Let $\ell$ be a prime number, and let $k$ be a finite field containing a primitive $\ell$-th root of unity. Hence, $\Char k \ne \ell$. 
Let $K$ be a \globalfieldv, and let $\curlyZ/k$ be a smooth complete curve with  function field $K$.
Let $C/K$, $C_D/K$, $J/K$, and $J_D/K$ be as defined in the case of number fields.
\end{definition}

Over a number field $K$, given any positive integer $g_0$, it does not seem to be known how to produce a (superelliptic) curve $C/K$ of genus $g(C) > g_0$ such that there are infinitely many twists with ($\lambda$-) Selmer rank $0$, or such that the Selmer rank for $C/K$ is $0$. We use Theorem \ref{thm:functionfield1} below to show that a function field analogue of this problem can be solved even for (infinitely many) function fields $K$ with arbitrarily large genus $g(K/k)$.

\medskip
\copytheorem{Theorem \ref{thm:functionfield1}}
Assume the hypotheses in \ref{def:ff-setup}.
Suppose that $f(x)$ is defined over the constant field $k$, irreducible over the finite field $k$, and has prime degree $p$.  Let $k'$ be the field extension of $k$ of degree $p\ne \Char k$. Let $L:=K\otimes k'$, and let  $\OL$ be the integral closure of $\OK$ in $L$.
\par
Let $A/k$ be the Jacobian variety of the normalization of the superelliptic curve $y^\ell =f(x)$ defined over the constant field $k$, so $A_K=J$.
If $\Sharp\, \Cl(\OL) \not\equiv 0 \mod \ell$, then $\dimp\ell\sellam(J,K)=0$.  Moreover, there is a set $\curlyD$ of prime ideals of $\OK$ with Dirichlet density $(p-1)/p$ such that whenever $D$ is a nonzero element of $\OK$ supported by $\curlyD$,
$$\dimp\ell\sellam(J_D,K)=0,\ \text{ and }\ %
                \Sharp\ C_D(K) \le \Sharp\ A(k).$$
\rm
\indent
The hypothesis in Theorem \ref{thm:functionfield1} is often satisfied.
Let $K$ be a function field of one variable with a rational divisor $v_\infty$, and let $\curlyZ/k$ be the smooth curve with function field $K$.
Then, by Proposition \ref{prop:dino}, $\Sharp\, \Cl(\OK) = \Sharp\, \Pico(\curlyZ)$. Let $L:=K\otimes k'$ for a finite extension $k'$ of $k$, and let $\OL$ be the integral closure of $\OK$ in $L$. Then $\Sharp\, \Cl(\OL) = \Sharp\, \Pico(\curlyZ_{k'})$. 
Therefore, for all but finitely many prime numbers $\ell$, we have $\Sharp\, \Cl(\OL) = \Sharp\, \Pico(\curlyZ_{k'})\not\equiv 0 \mod \ell$.
Thus, we find examples of superelliptic curves of arbitrarily large genus which satisfy the conditions in Theorem \ref{thm:functionfield1}.
\par
For the curves $C$ considered in Theorem \ref{thm:functionfield1}, there are infinitely many $D$'s such that $\Sharp\,C_D(K)$ is bounded. 
In \cite{schoen:1990}, Schoen considered hyperelliptic curves defined over  certain geometric fields, and showed that for these curves, the number of rational points of their quadratic twists can be arbitrarily large.

\par
We conclude this introduction by providing the reader with a road map for Schaefer's description for the $\lambda$-Selmer group, and for the proof of Theorem \ref{MainTheorem1}.
Recall that in this case, $f(x)$ is irreducible over $K:=\ratn(\zetal)$.  Let $L$ be a field isomorphic to $K[x]/(f)$.  Using Schaefer's method, we have the first two rows of the commutative diagram (\ref{diag:together}).
\begin{figure}[b]
\begin{equation}\label{diag:together}
\xymatrix@R=\baselineskip@C=50pt{
        \WMell \ar[r]^\delta \ar[d]^{\kappav} 
                & \HonE(K,J[\lambda])_{S_J} \ar[r]^{\theta} \ar[d]^{\res_v}  
                        & L(S_J,\ell) \ar[r]^{\incl} \ar[d]^{\res_v}\ar[dr]^{\Psi_J} 
                                & \thegroupell{L}\ar[d]^{\normAX[L]}\\
        \WMell[\Kv] \ar[r]^{\delta_v} \ar@{.>}[d]^{\curlyh_v^D} 
                & \HonE(K_v,J[\lambda]) \ar[r]^{\theta_v} \ar[d]^{\id^D} 
                        & \thegroupell{L_v} \ar@{=}[d] 
                                & \thegroupell{K}\ar@{=}[d]\\                        
        \WMDell[\Kv] \ar[r]^{\delta_v^D} 
                & \HonE(K_v,J_D[\lambda]) \ar[r]^{\theta_v^D} 
                        & \thegroupell{L_v}  
                                & \thegroupell{K}\\
        \WMDell \ar[r]^{\delta^D} \ar[u]_{\kappav} 
                & \HonE(K,J_D[\lambda])_{S_D} \ar[r]^{\theta^D} \ar[u]_{\res_v}  
                                & L(S_D,\ell) \ar[r]^{\incl} \ar[u]_{\res_v}\ar[ur]^{\Psi_D} 
                                        & \thegroupell{L}\ar[u]_{\normAX[L]}\diadot
        }
\end{equation}
\end{figure}
\ \\
\indent
Recall from \ref{def:CDJD} that a twist $J_D/K$ is the Jacobian variety of $C_D/K$ for some $D\in \OK\minuS\set{0}$, and the curve $C_D$ is given by $y^\ell = f_D(x):=D^p\, f(x/D)$. Since $L \cong K[x]/(f_D)$, as in the case of $J/K$, we can construct a map $\theta^D : \HonE(K,J_D[\lambda]) \to \thegroupell{L}$.
When we consider the $\lambda$-Selmer groups of both $J$ and $J_D$, we establish that the first two rows and the last two rows of \togeth\ are commutative. It is noteworthy that the targets of $\theta$ and $\theta^D$ are both $\thegroupell{L}$, and it can be understood as a consequence of the fact that the two group schemes $J[\lambda]$ are $J_D[\lambda]$ are isomorphic to each other.
\par
Two key points we shall prove in Section \ref{sec:def-prop} are the following: 
First, if $D$ is an $\ell$-th power in $K_v$ for $v$, then there is a map $\curlyh_v^D : J(K_v) \to J_D(K_v)$
such that the diagram \togeth\ commutes (see Proposition \ref{prop:key2}).  
It is clear that there is an isomorphism: $J(K_v) \to J_D(K_v)$, but it is not so obvious, unless one knows the horizontal maps, that it commutes with the identiy map on $\thegroupell{L_v}$.
Secondly, if $D$ is divisible only by prime ideals $\primep$ of $\OK$ such that $\primep\OL$ is prime, then $\ker \Psi_J =\ker\Psi_D$ in the diagram \togeth\ (see Theorem \ref{Sec:proof:prop2}).
These two results immediately imply that $\theta^D( \sellam(J_D,K) ) \subset \theta( \sellam(J,K) )$ (see the proof of Theorem \ref{MainTheorem1}).
By posing more conditions on $D$, we can prove
$\dimp\ell \sellam(J_D,K) = \dimp\ell \sellam(J,K)$. 
To exhibit the existence of (enough) such $D$'s, we use the Cebotarev density theorem and the general reciprocity laws.
\paragraph{\bf Notation}
 \begin{definition}\label{rel-infty}
\textit{A global field} is either a number field or a function field of one variable over a finite field $k$.
Let $K$ be a global field. We denote by $\MK$ the set of all places of $K$.
Suppose that $K$ is a function field of one variable over a finite field $k$. We choose a finite set of places of $K$, which will be denoted by $\MKinf$, and consider the Dedekind domain 
$$ R:=\set{ \al \in K^* : \abs{\al}_v \le 1,\text{ for all } v \in \MK \backslash \MKinf}.$$
We shall denote $\MK \backslash \MKinf$ by $\MKfinite$.  
\par
In this paper, by a global field, we shall mean a number field with $\OK$, the integral closure of $\zz$, or a function field of one variable over a finite field $k$ with the Dedekind domain $\OK$ determined by a choice of  $\MKinf$.
For both cases, $\OK$ is called the ring of integers of $K$. 
Throughout the paper, if $L$ is a finite separable extension of $K$, let $\OL$ denote the integral closure of $\OK$ in $L$.
  When $K$ is a function field of one variable, our choice of $\MKinf$ determines $\OK$ and, hence, $\OL$ and $M_L^\infty$.
For a function field $K$ of one variable with a rational divisor $v_\infty$ we choose $\MKinf:=\set{v_\infty}$.
\end{definition}

\begin{definition}\label{def:GK-GF}
Throughout the paper, given a field $K$, let $\Kbar$ denote the algebraic closure of $K$, and  $\Ksep$, the (algebraic) separable closure of $K$.  Let $K$ be a global field.
We denote by $\GGK$ the absolute Galois group $\GK$. 
For each $\primeq\in M_K$, let $\Kq$ denote the completion of $K$ at $\primeq$. 
We fix the algebraic closures $\Kbar$ and $\Kbar_\primeq$, and fix an embedding $\kappa_\primeq : \Ksep \injects \Kqsep$.  
\end{definition}

In this paper, \textit{a variety $J$ over an arbitrary field $K$} is a separated scheme of finite type over $K$ such that $J_{\Kbar}$ is integral (that is, reduced and irreducible).
Let $J/K$ be a  variety.
For a field extension  $F$ of $K$, \textit{the $F$-points of $J$} is the set  $J(F)$.
We denote by $\overline J$ the variety $J_{\Kbar}$, and by $J\sep$ the variety $J_{\Ksep}$.

\section{Schaefer's description of the $\lambda$-Selmer group}\label{subsub:Selmer}
In this section, we  introduce Schaefer's method for computing Selmer groups.  This method is introduced in \cite{schaefer:1998} for number fields; however, his proofs carry over to the case of global fields.

\begin{definition}\label{def:l-des}
Let $\ell$ be a prime number, and let $K$ be a global field of characteristic $\ne \ell$ such that $K$ contains a primitive $\ell$-th root  of unity $\zetal$. 
Let $C/K$ denote the normalization of a \super\ given by $y^\ell=f(x)$ with $d:=\deg(f)$.
Let $J/K$ be the Jacobian variety of $C$.
Let $\zetal$ also denote the $K$-automorphism on $J$  defined in \cite{schaefer:1998}, Sec 3, and $\lambda:=1-\zetal \in \End(J)$, which we also denote by $[\lambda]$.
Recall that $f(x)$ has distinct roots, and let $\set{z_i\in K\sep : i=1,\dots,d}$ be the roots of $f(x)$.  
Let 
$P_i:=[(z_i,0) - (\infty)]  \in J(K\sep)$.  By \cite{schaefer:1998}, 
Proposition 3.2, the set $\set{P_1,\dots,P_d}$ generates the $\finiteell$-vector space $J[\lambda](K\sep)$ of dimension $d-1$. 
We choose a $\GGK$-set $X$ as small as possible such that $X$ spans $J[\lambda](K\sep)$:
If $f(x)$ \noroot, then we set $d':=d$, and if $f(x)$ has a $K$-rational root $z_d$, then we set $d':=d-1$. Then we define
$X := \set{P_1,\dots,P_{d'}}$.  The $\Ksep$-algebra of all functions from $X$ to $\Ksep$ is denoted by $\Abar$.  As described in \cite{schaefer:1998}, Sec 2, $\Abar$ has $\GGK$-action, and the $\GGK$-invariants are denoted by $A_X$.
\end{definition}

\begin{definition}\label{def:Sel}
From the Galois cohomology for the sequence
$0 \to J[\lambda](\Ksep) \to J(\Ksep) \overset{\lambda}{\to} J(\Ksep) \to 0$, 
we have  the following injective map: 
\begin{equation}\label{prop:descent:eq}
\delta : \WMell \To \HonE(K,J[\lambda]).
\end{equation}
\par
For each $v\in M_K$,  we denote by $\res_v$ the map: $\HonE(K,J[\lambda])\to\HonE(K_v,J[\lambda])$, and  by $\delta_v$ the coboundary map: $J(K_v) \to \HonE(K_v,J[\lambda])$ defined for $K_v$.
\textit{The $\lambda$-Selmer group of $J/K$} is 
\begin{equation}\label{eq:def-Sel}
\Sel{\lambda}(J,K):=\set{\xi \in \HonE(K,J[\lambda]) : \res_v(\xi)\in \Img \delta_v, \text{ for all } v\in M_K}.
\end{equation}
\end{definition}
By definition, \textit{ the $\lambda$-part of the Tate-Shafarevich group of $J/K$}, denoted by $\sha(J/K)[\lambda]$, is canonically isomorphic to $\Sel{\lambda}(J,K) / \Img \delta$. 
Let $S$ be a subset of $\MK$ containing all archimedean places of $K$.
The subgroup of $\HonE(K,J[\lambda])$ unramified outside the set $S$ is denoted by
 by $\HonE(K,J[\lambda])_S$.
Then the following result is standard:
Let $S$ be a subset of $\MK$ containing all \arch\ places, the places above $\deg(\lambda)$, and the places of bad reduction of $J/K$. 
Then,
\begin{equation}\label{prop:Sel}
\Sel{\lambda}(J,K)=\set{ \xi \in \HonE(K,J[\lambda])_S : \res_v(\xi) \in \Img \delta_v \text{ for all } v \in S}.
\end{equation}
Moreover, $\Sel{\lambda}(J,K)$ is finite, and contains $\delta( J(K)/\lambda J(K) )$.

\par
Suppose that $f(x)$ \noroot.
Then $\dimp\ell \ellthroots = d$ since $\Sharp X=d$, and that $\dimp\ell J[\lambda](K\sep)=d-1$.
By \cite{schaefer:1998}, Proposition 3.4,
the following is a split short exact sequence of $\GGK$-modules, which is due to the fact that $\ell \nmid d$:
 \begin{equation*}
 0 \To J[\lambda](K\sep) \overset{w}{\To} \ellthroots \overset{\Norm_{\Abar/\Kbar}}{\To} \ellthroots[K\sep] \To 0
 \end{equation*}
 where $w$ is the map defined in \cite{schaefer:1998}, page 452.
 Since this exact sequence splits, the following natural map on the cohomology groups is injective:
\begin{equation}\label{thm:injective}
\til w : \HonE(K,J[\lambda]) \To \HonE(K,\ellthroots).
\end{equation}
\par
Suppose that $f(x)$ has a $K$-rational root, say $z_d$.  
We choose $X:=\set{P_1,\dots,P_{d-1}}$ as in \ref{def:l-des}.
Then it is clear that $w : J[\lambda](K\sep) \to \ellthroots$ is an isomorphism of $\GGK$-modules since $\dimp\ell J[\lambda](K\sep) = \dimp\ell \ellthroots=d-1$.
Therefore,  the following natural map is an isomorphism:  
\begin{equation}\label{eq:small-X}
\til w : \HonE(K,J[\lambda]) \To \HonE(K,\ellthroots).
\end{equation}

\begin{definition}\label{def:Kummer}
Consider the Kummer sequence for $\Abarstar$ with respect to the homomorphism $\al \mapsto \al^\ell$. By  Hilbert Theorem 90, we have a natural isomorphism:
$\Phi : \HonE(K,\ellthroots)\to \thegroupell{A_X}$ and, hence, have an injective homomorphism
\begin{equation}\label{TheInjection}
\Phi\circ \til w : \HonE(K,J[\lambda]) \To \thegroupell{A_X}.
\end{equation}
\end{definition}

\par
Let $S$ be a subset of $\MK$ containing $\MKinf$, places above $\ell$, and places of bad reduction of $J/K$.
If $f(x)$ has a $K$-rational root, then  
$\HonE(K,J[\lambda])_S \cong \AX(S,\ell)$
where $\AX(S,\ell)$ is the subgroup of $\thegroupell{\AX}$ defined in \cite{schaefer:1998}, Sec 2.
If $f(x)$ does not have a $K$-rational root, then, \cite{schaefer:1998}, Proposition 3.4 gives us the following isomorphism:
\begin{equation*}
\HonE(K,J[\lambda])_S \cong \ker\big( \normAX : \AX(S,\ell) \to K(S,\ell) \big).
\end{equation*}

\begin{definition}\label{def:kappa}
For each $v\in M_K$, 
we take the image of $X$ in $J[\lambda](\Kvsep)$ to be a $\GGKv$-stable spanning subset, and denote it by $\Xv$.
Then we have an induced map
\begin{equation}\label{subsub:bij}
 \Abar \To \overline{A}_{\Xv}
\end{equation} 
obtained by pulling back the natural map $\Xv \to X$, and 
it is denoted by $\til \kappav^*$. Moreover, we have the maps $\delta_v$, $\til w_v$, and $\Phi_v$ as in the case of $K$.
\end{definition}

\par
Let $S$ be a set of the places containing $\MKinf$, the places above $\ell$, and the places of bad reduction of $J/K$.
For each $v\in M_K$, 
the following diagram is commutative:
\begin{equation*}
\xymatrix@C=60pt{
        \WMell \ar[r]^\delta\ar[d]_{\kappa_v} & \HonE(K,J[\lambda])_S \ar[d]_{\res_v} \ar[r]^{\Phi\circ \til w} 
                                                                                        & \AX(S,\ell)\ar[d]_{\til\kappa_v^*}\\
        \WMell[K_v] \ar[r]^{\delta_v} 
                                                        & \HonE(K_v,J[\lambda]) \ar[r]^{\Phi_v\circ \til w_v} 
                                                                 & \thegroupell{A_{X_v}} \diadot    
                                                                 }
\end{equation*}

\par
Suppose that $f(x)$ \noroot. Then 
$$\sellam(J,K)\cong
        \set{\al \in \AX(S,\ell) : \normAX(\al) = 1,\ \til\kappa_v^*(\al)\in\Img \Phi_v\circ \til w_v \circ \delta_v \text{ for all }
                v \in S }.$$
Suppose that $f(x)$ has a $K$-rational root. 
Then, we have a simpler description:
$$\sellam(J,K)\cong
        \set{\al \in \AX(S,\ell) :  \til\kappa_v^*(\al)\in\Img \Phi_v\circ \til w_v \circ \delta_v \text{ for all }
                v \in S }.$$

Let $\Phi $ and $\til w $ be the maps defined in (\ref{TheInjection}).  Let $\delta $ be 
the coboundary map defined in (\ref{prop:descent:eq}). Then, we have an injective homomorphism
\begin{equation}\label{eq:descr-cob}
\Phi\circ\til w \circ \delta : \WMell[K] \To \thegroupell{\AX},
\end{equation}
and there is a useful description of this map which given below.
\par
Write $K\sep(C\sep)$ as the field of fractions of $K\sep[x,y]/(y^\ell - f(x))$.
Then $x-z_i$ is an element of this function field, which we will denote by $f_{P_i}(x,y)$; 
\begin{equation}\label{eq:fPi}
f_{P_i}(x,y):=x-z_i.
\end{equation}
Note that $f_{P_i}$ can be thought as a function on divisors of $C$. 
If $E$ is a divisor in $\Div(C)$ \textit{avoiding} the set $X$ (see \cite{schaefer:1998}, page 450), then we denote by $f_\bullet(E)$ the function: $X \to \Kbar^*$ given by $P \mapsto f_P(E)$ for all $P\in X$. 
Let $Y:=\set{P_1,\dots,P_d}$.  If $E$ is a divisor in $\Div(C)$ avoiding $Y$, then, as long as  there is no confusion, we also denote by 
$f_\bullet(E)$ the extended function: $Y \to \Kbar^*$ given by $P \mapsto f_P(E)$ for all $P\in Y$.
If $E:=\sum (R)$ is a divisor in $\Div(C)$ avoiding $X$ such that $K(R)/K$ is separable, then $\fdot(E) \in \Astar$, and if  $E$ is a divisor in $\Div(C)$ avoiding $Y$ such that $K(R)/K$ is separable, then $\fdot(E) \in \AX[Y]^{\ *}$.

\begin{definition}\label{def:fdotE}
Let us define a group homomorphism $\Pico(C) \to \thegroupell{\AX}$ as follows:
It is well-known that for each $[D] \in \Pico(C)$, we can choose a divisor $E=\sum (R)$ avoiding any given finite set of $C$ such that $K(R)/K$ is separable and $[E]=[D]$.  Then we define 
a group homomorphism given by
\begin{equation}\label{eq:def-coboundary}
[ D ] \mapsto \cls(f_\bullet(E)) \in \thegroupell{\AX},
\end{equation}
and denote it by $\delbar$.
\end{definition}
\par
The following corollary follows from \cite{schaefer:1998}, Theorem 2.3 and Proposition 3.3.  The proposition essentially shows that $[(\infty)]=[D]$ for some divisor $D$ avoiding $Y$ such that $f_\bullet(D) \in (\Astar[Y])^\ell$:
\begin{lemma}\label{cor:coboundary}
Let $\delbar$ be the map defined in \ref{def:fdotE} for $X$.
Then $\Phi  \circ \til w  \circ \delta = \delbar$.
\par
Suppose that $f(x)$ has a $K$-rational root $z_d$, so we choose  $X:=\set{P_1,\dots,P_{d-1}}\subset J(K\sep)$.
If $E$ is a nonzero divisor in $\Divo(C)$ such that $E=(Q)-(\deg_K(Q))(\infty)$ for some $Q$ avoiding $X$, and such that $K(Q)/K$ is separable, then 
$[E]$ is mapped to $\cls(f_\bullet(Q))\in \thegroupell{\AX}$ under $\delbar$.  In particular, $[(z_d,0)-(\infty)]\mapsto \cls( f_\bullet(z_d,0) ) \in \thegroupell{\AX}$.
\end{lemma} 
\begin{proof} The proof is left to the reader.
\end{proof}

\section{Two Key Propositions}\label{sec:def-prop}

Let us recall the diagram \togeth. In this section, we prove the commutativity of this diagram in a slightly more general context. 
Let $\ell$ be a prime number, and let $K$ be a global field of \chr\ $\ne \ell$, containing a primitive $\ell$-th root $\zetal$ of unity.  Let $C/K$ be the normalization of a \super\ given by $y^\ell=f(x)$, and let $d:=\deg(f)$.
\begin{definition}\label{def:XD}
For each $\primeq\in \MK$, we fix an embedding $\kappa_\primeq : \Kbar \to \Kbar_\primeq$.
Let $X$ and $X_\primeq$ for each $\primeq\in\MK$ be the $\GGK$-stable spanning sets of $J[\lambda](K\sep)$ and $J[\lambda](\Kq,\sep)$, respectively, defined in \ref{def:l-des} and \ref{def:kappa}.
For each $\primeq\in\MK$, let $\kappa_\primeq$  and $\til\kappa_\primeq^*$ be the maps 
defined in Section \ref{subsub:Selmer};
\begin{equation*}
\kappa_\primeq : X \To X_\primeq,\quad 
        \til\kappa_\primeq^* : \Abar \To \Abar[X_\primeq].
\end{equation*}
Let $P_i$ for $i=1,\dots,d'$ be the points in $X$.
Note that $f(x)=(x-z_1)\cdots(x-z_d)$, and for $D \in K^*$, $D^d f(x/D) = (x-z_1\, D)\cdots(x-z_d\, D)$.
For $D \in K^*$ and $\primeq \in \MK$, we define
\begin{align*}
X^D &:=\set{ [(z_i\,D,0)-(\infty)]\in J_D(K\sep)  : i=1,\dots,d'};\\
\Xq^D &:= \set{ [(\kappaq(z_i\,D),0)-(\infty)] \in J_D(\Kq,\sep) : i=1,\dots,d'};\\
P_i^D &:=[(z_i\,D,0) - (\infty)] \in J_D(K\sep),
\end{align*}
and let $\Phi^D$, $\Phi^D_\primeq$, $\til w^D$,  $\til w^D_\primeq$, $\delta_\primeq$, and $\delta_\primeq^D$  be the maps for the  twist $J_D$ as in Section \ref{subsub:Selmer}.
Then we have the following sequences of injective maps:
\begin{equation}\label{eq:d-phi}
\xymatrix@R=\baselineskip@C=40pt{
        \WMDell\ar[r]^{\delta^D}\ar[d]^{\kappa_\primeq} & \HonE(K,J_D[\lambda])\ar[r]^{\til w^D} \ar[d]^{\res_\primeq}
                                                                & \HonE(K,\ellthroots[{\Abar[\XD]}])\ar[r]^{\Phi^D}\ar[d]^{\res_\primeq} 
									& \thegroupell{A_{X^D}};\ar[d]^{\til\kappa_\primeq^*} \\
        \WMDell[\Kq]\ar[r]^{\delta_\primeq^D} & \HonE(\Kq,J_D[\lambda])\ar[r]^{ \til w^D_\primeq}                                              
                                                                   & \HonE(\Kq,\ellthroots[{\Abar[\XDq]}])\ar[r]^{\Phi^D_\primeq} 
                                                                                                        & \thegroupell{A_{X_\primeq^D}}. 
        }                                                                                                        
\end{equation}        
\end{definition}

\par
 Let  $Z(J):=\set{T_i \in X: i=1,\dots,s}$ be a set of representatives of $\GGK$-orbits in $X$, and $Z(J_D):=\set{P_j^D \in X^D : P_j\in Z(J)}$ which is said to be \textit{compatible with} $Z(J)$.  Let $L_i$ be the subfield of $K\sep$, generated by $T_i$ for $i=1,\dots,s$.
Since the action of $\GGK$ is made through the divisors representing the points in $X$, 
the subfield $K(P_j^D)$ of $K\sep$ is $K(z_j)$, and $Z(J_D)$ is a set of representatives of $\GGK$-orbits in $X^D$.  
For $\primeq \in \MK$, let $Z(J,\primeq)$ be a set of representatives of $\GGKq$-orbits such that $Z(J) \subset Z(J,\primeq)$. Then, we choose $Z(J_D,\primeq)$ compatible with $Z(J,\primeq)$; hence, $Z(J_D) \subset Z(J_D,\primeq)$.

\begin{definition}\label{def:C-Cq}
Note that if $\al\in \AX[X^D]$, then the evaluation of $\al$ at $P_j^D\in Z(J_D)$ is contained in $K(z_j)=L_i$ for some $i$, i.e., $\al(P_j^D) \in L_i$.
For $D\in K^*$ and $\primeq \in \MK$, let $\Psi_{Z(J_D)}$ and $\Psi_{Z(J_D,\primeq)}$ be the evaluation maps of 
 $\thegroupell{\AX[\XD]}$ and $\thegroupell{\AX[X_\primeq^D]}$ at $Z(J_D)$ and $Z(J_D,\primeq)$, respectively, and define
\begin{equation*}
\begin{aligned}\label{eq:Z-ZD}
\curlyC &:=\Psi_{Z(J_D)}\big( \thegroupell{A_{X^D}} \big) 
        = \prod_{ P \in Z(J_D) } \thegroupell{ K(P) }
        = \prod_{i=1}^s \thegroupell{L_i};\notag\\
\Cq &:= \Psi_{Z(J_D,\primeq)}\big( \thegroupell{A_{X_\primeq^D}} \big) 
        = \prod_{P\in Z(J,\primeq)} \thegroupell{\Kq(P)}.
\end{aligned}
\end{equation*}
\end{definition}
Since we choose $Z(J_D,\primeq)$ compatible with $Z(J,\primeq)$, the group $\Cq$ is defined not depending on $D$ as  $\curlyC$.

\par
Recall the maps in (\ref{eq:d-phi}).
For each $D\in K^*$, we denote by $\theta$ and $\theta^D$ the injective maps $\Psi_{Z(J)} \circ \Phi \circ \til w$ and
$\Psi_{Z(J_D)} \circ \Phi^D \circ \til w^D$, respectively.
For each $D\in K^*$ and $\primeq\in\MK$, we denote, respectively, by $\theta_\primeq$ and $\theta_\primeq^D$ the injective maps
\begin{gather}
\Psi_{Z(J,\primeq)}\circ\Phi_\primeq\circ \til w_\primeq : \HonE( K_\primeq,J[\lambda]) \To \Cq;\label{eq:w-Phi}\\
\Psi_{Z(J_D,\primeq)}\circ\Phi^D_\primeq\circ \til w^D_\primeq : \HonE( K_\primeq,J_D[\lambda]) \To \Cq.\label{eq:w-Phi-D}
\end{gather}

\begin{remark}\label{rm:isomorphism}
It is noteworthy to observe that the targets  of the following maps are defined depending only on $C/K$:
\begin{equation*}\label{eq:C-indep}
\theta^D : \HonE(K,J_D[\lambda]) \To \curlyC,\quad 
\theta_\primeq^D : \HonE(\Kq,J_D[\lambda]) \To \Cq.
\end{equation*}
\par
Let $S$ be a subset of $\MK$ containing $\MKinf$, the places above $\ell$, and the places of bad reduction of  $J/K$ and $J_D/K$. If $f(x)$ \noroot, then 
\begin{equation*}
\theta^D(\HonE(K,J_D[\lambda])_S) 
        = \ker\big( \pprod_i\normAX[L_i] : \pprod_i L_i(S,\ell) \to K(S,\ell) \big) 
         = \theta(\HonE(K,J[\lambda])_S).
\end{equation*}
If $f(x)$ has a $K$-rational root, then 
\begin{equation*}
\theta^D\big( \HonE(K,J_D[\lambda])_S \big) 
        = \pprod_{i=1}^s L_i(S,\ell) = \theta\big( \HonE(K,J[\lambda])_S \big) .
\end{equation*}    
\end{remark}
\par
Recall the diagram \togeth\ from the introduction.  In a slightly more general context, we constructed above all the horizontal maps in \togeth.  Proposition \ref{prop:key1} below is one of the key propositions which will establish the commutativity of the diagram formed by  the restriction maps in the second and the third columns of \togeth. The proof is left to the reader.
\begin{proposition}\label{prop:key1}
For each $\primeq\in M_K$, there is a natural  map $\curlyC \to \Cq$ such that the following diagram commutes for all $D\in K^*$:
\begin{equation*}
\xymatrix@C=50pt{
\HonE(K,J_D[\lambda]) \ar[d]^{\res_\primeq} \ar[r]^{\til w^D} 
        & \HonE(K,\ellthrootS[{\Abar[X^D]}])\ar[d]^{\res_\primeq} \ar[r]^{\Phi^D} 
                & \thegroupell{A_{X^D}} \ar[r]^{\Psi_{Z(J_D)}} \ar[d]_{\til\kappa_\primeq^*} 
                        & \curlyC\ar[d]_{\res_\primeq}\\
\HonE(\Kq,J_D[\lambda])  \ar[r]^{\til w_\primeq^D} 
        & \HonE(\Kq,\ellthrootS[{\Abar[X_\primeq^D]}]) \ar[r]^{\Phi_\primeq^D}          
                & \thegroupell{A_{X_\primeq^D}} \ar[r]^{\Psi_{Z(J_D,\primeq)}}  
                        & \Cq\diadot
        }
\end{equation*}
\end{proposition}
\par
Proposition \ref{prop:key2} below completes the proof of the commutativity of the diagram \togeth.
\begin{proposition}\label{prop:key2}
Let $\primeq$ be a place in $\MK$. For all 
 nonzero elements $D$ of $\OK$  such that $D \in \ellthpower{\Kq}$,
 there is an isomorphism $\curlyh_D : J(\Kq) \to J_D(\Kq)$ such that the following is a commutative diagram:
\begin{equation}\label{diag:D1D2}
\xymatrix@C=90pt{
        \WME{J_{ D}}{ K_\primeq} \ar[r]^{\theta_\primeq^{  D }\delta_\primeq^{  D}}
                        &  \Cq\\
        \WME{J_{ }}{ K_\primeq} \ar[r]^{\theta_\primeq^{ }\delta_\primeq^{ }}\ar[u]^{\curlyh_D} 
                                                &  \Cq\ar[u]^{\text{\rm id}}\diadot
        }
\end{equation}
In particular, $\Img \theta_\primeq^D \delta_\primeq^D=\Img \theta_\primeq \delta_\primeq\ $ for all nonzero $D \in \OK$ and $\primeq\in\MK$ such that $D \in \ellthpower{\Kq}$.
\end{proposition}
\begin{proof}
Let $D$ be a nonzero element of $\OK$ such that $D \in \ellthpower{\Kq}$.
Then, we have an isomorphism $\CDKq\to \CKq$ given by $(x,y) \mapsto (x/D, y/\sqrt[\ell]{D^d})$, and this isomorphism induces an isomorphism $\curlyh_D : J(\Kq) \to J_D(\Kq)$ by pulling back on the divisors.
\par
Recall $Z(J_D,\primeq) =\set{P^D:=[(zD,0)-(\infty)] : [(z,0)-(\infty)] \in Z(J,\primeq) }$.
Let $E:=\sum n_j(R_j)$ be a divisor in $\Divo(\CKq,\sep)$ avoiding $X$ such that $K(R_j)/K$ are separable, and write 
$R_j=(x_j,y_j)$. Then, 
\begin{align*}
 f_\bullet\big( \curlyh_D(E) \big) (P^D) 
        &=\prod_j\big( f_{P^D}(x_j\,D, y_j\, \sqell{D^d})\big)^{n_j} 
        =\prod_j (x_j\, D - z\,D)^{n_j} \\
        &=D^{\sum n_j}\,\prod_j\ (x_j - z)^{n_j}
        = \prod_j\ (x_j - z)^{n_j};\\
 f_\bullet(E)\  (P) 
        &= \prod_j\ \big( f_{P} (x_j,y_j)\big)^{n_j}
                = \prod_j\ (x_j - z)^{n_j}.
\end{align*}
This proves the commutativity of the diagram (\ref{diag:D1D2}).\end{proof}

\begin{remark}\label{rem:idD}      
The diagram in Proposition \ref{prop:key2} can be put into the following diagram:
\begin{equation*}
\xymatrix@C=65pt{
        \WMDell[ K_\primeq] \ar[r]^{\delta_\primeq^{  D}}
                & \HonE(\Kq,J_D[\lambda]) \ar[r]^{\Phi_\primeq^D \circ \til w_\primeq^D} 
                        & \thegroupell{\AX[\Xq^D]} \ar[r]^{\Psi_{Z(J_D,\primeq)}}
                                &  \Cq\\
        \WMell[K_\primeq] \ar[r]^{\delta_\primeq^{ }}\ar[u]^{\curlyh_D} 
                & \HonE(\Kq,J[\lambda]) \ar[r]^{\Phi_\primeq \circ \til w_\primeq}\ar[u]^{\til\curlyh} 
                        & \thegroupell{\AX[\Xq]} \ar[r]^{\Psi_{Z(J_D,\primeq)}} \ar[u]^{\id_D}
                                                &  \Cq\ar[u]^{\text{\rm id}} \diadot
        }
\end{equation*}
where $\til\curlyh$ and $\id_D$ are  isomorphisms induced from the natural isomorphism: $J[\lambda] \to J_D[\lambda]$.

All the vertical maps in the above diagram, except the first one, exist whether or not $D$ is an $\ell$-th power in $\Kq$, and the entire diagram commutes when $D \in \ellthpower{\Kq}$.
\end{remark}

\begin{definition}\label{Sel2}
Recall the maps $\theta^D$ and $\theta_\primeq^D$ for $\primeq \in \MK$ and $D \in \OK$ defined in \ref{def:C-Cq}.
Let $S$ be a subset of $\MK$ containing $\MKinf$, the places above $\ell$, and the places of bad reduction of $J_D/K$.
Then, 
$\sellam(J_D,K)$ is described as a subgroup of $\curlyC$ as follows:
\begin{equation*}
\theta^D(\sellam(J_D,K)) = 
        \set{ \al \in \theta^D\big( \HonE(K,J_D[\lambda])_S \big) : \res_\primeq(\al) \in \Img \theta_\primeq^D\delta_\primeq^D
                \text{ for all } \primeq \in S }.
\end{equation*}
\end{definition}

\section{The Jacobian Varieties Without $\lambda$-torsion Points}\label{Sec:2-descent:twist}

In this section, we prove Theorem \ref{MainTheorem1} and its analogue for a function field.
Let $\ell$ be a prime number, and let $K$ be a \globalfield.  
Let $C/K$ be the normalization of a \super\ $y^\ell=f(x)$.  Let $\Delta_f$ be the discriminant of $f(x)$
where $f(x)$ is irreducible over $K$ with prime degree $p$.
Let $J/K$ be the Jacobian variety of $C$.  We keep all the notation and definitions introduced in Section \ref{sec:def-prop}.

\begin{definition}\label{SJ-SD}
Let $S_J$ denote the subset of $\MK$ containing $\MKinf$, the places dividing $\ell\Delta_f$.
For each nonzero element $D$ of $\OK$, let $\SD$ denote the set $\SE\cup\set{\primep \in \MK^0 : \primep \mid D\OK}$. 
Then, the sets $S_J$ and $S_D$ contain the set of places of bad reduction of $J/K$ and $J_D/K$, respectively.
\end{definition}

\begin{definition}\label{def:context1}
Recall that  $X:=\set{P_1,\dots,P_d}$ is the subset of $J[\lambda](K\sep)$.
 Since $X$ forms one orbit,  $Z(J)$ contains a single point, and 
we choose $Z(J):=\set{P_1}$ as a representative.  Let $L$ denote the finite separable field extension $L_1:=K(P_1)$ of degree $p$. Then $\curlyC:=\thegroupell{L}$ as defined in \ref{def:C-Cq}.
\end{definition}

\begin{theorem}\label{Sec:proof:prop2}
Let $D$ be a nonzero element of $\OK$, and suppose that $D\OK$ is supported by the set of  prime ideals $\primeq$ of $\OK$ such that either $\primeq\OL$ is prime in $\OL$ or $\primeq\OL=\primep^p$ for some prime ideal $\primep$ of $\OL$.
Then,  
$\ker\big( N_{L/ K} : L(S_J,\ell) \to  K(S_J,\ell) \big)=\ker\big( N_{L/ K} : L(S_D,\ell) \to  K(S_D,\ell) \big)$.
Hence,    $\theta\big( \HonE( K,J[\lambda])_{S_J} \big) = \theta^D\big( \HonE( K,J_D[\lambda])_{S_D} \big)\text{ as subgroups of  } \thegroupell{L}$.
\end{theorem}
\begin{proof}
Since $S_J\subset S_D$, it is clear that 
\begin{equation*}
\ker\big( N_{L/ K} : L(S_J,\ell) \to  K(S_J,\ell) \big)
        \subset \ker\big( N_{L/ K} : L(S_D,\ell) \to  K(S_D,\ell) \big).
\end{equation*}
\par
Let $\al$ be a nonzero element of $\OL$ such that $\al\OL = \idealA\idealB^\ell$ where $\idealA$ is an ideal \textit{supported by} $S_D$ and such that $N_{L/ K}(\al)=\beta^\ell$ for some $\beta\in \OK$.  Let $\primeq$ be a prime ideal of $\OK$ dividing $D\OK$  such that $\primeq\OL$ is $\primep$ or $\primep^p$ for some prime ideal $\primep$ of $\OL$.
Let $n:=\ord_\primep(\al\OL)$, and let $m$ be the residue degree  of $\primep$ with respect to $\primeq$.  Since $m\not\equiv 0 \mod \ell$, and there is only one prime ideal of $\OL$ lying over $\primeq$, it follows that 
$
\ord_\primeq(\beta^\ell\OK) 
        = \ord_\primeq(N_{L/ K}(\al\OL))
         =nm.
$
On the other hand, $\ord_\primeq(\beta^\ell\OK)\equiv 0 \mod \ell$ and, hence, $n\equiv 0 \mod \ell$. 
Therefore,  $\ord_\primep(\al\OL)\equiv 0 \mod \ell$ for all $\primep$ dividing $D\OL$, and  $\al\OL = \idealA' (\idealB')^\ell$ for some ideal $\idealA'$ supported by $S_J$ , i.e., $\cls(\al)\in L(S_J,\ell)$.\end{proof}
\par
The Legendre symbol is used throughout Sections \ref{Sec:2-descent:twist} and \ref{sec:full-torsion}. Definitions and results on Legendre symbols required for our main results are introduced and proved in Section \ref{sec:rec}.
In \ref{def:power-symbol}, we extend the definition of the symbol for prime ideals dividing $\ell\OK$ and archimedean places, so that given $\al\in\OK$ and $\primep \in M_K$,  $\legendre{\al}{\primep}=1$ implies $\al\in(K_\primep^*)^\ell$.

\subsection{The case of number fields}\label{case:number}
For Lemma \ref{lem:L-K}, Proposition \ref{Sec:proof:prop1}, and Lemma \ref{subsec:proof:lem2}, we assume a slightly more general context:
Let $K/\ratn$ be the $\ell$-th cyclotomic field extension of $\ratn$ where $\ell$ is a regular prime number. Let $L$ be a field extension of $K$ such that $p:=[L:K]$ is a prime number not equal to $\ell$.

\begin{definition} \label{Sec:proof:def3}
Let $W$ be a finite subset of $\OL\minuS\set{0}$.  We denote by $\curlyD_W$ the set of prime numbers $q\in\zz$ not dividing $\ell$ which satisfy the following properties: for all prime ideals $\primeq$ of $\OK$ dividing $q$, 
\begin{enumerate}
\item\label{D:property1} The ideal $\primeq\OL$ is prime;
\item\label{D:property2}  For all $\al\in W$, $\al\not\equiv 0 \mod \primeq\OL$, and 
        $\legendre{\al}{\primeq\OL}=1$.
\end{enumerate}
\end{definition}

\begin{lemma}\label{lem:L-K}
Let $\primeq$ be a prime ideal of $\OK$ not dividing $\ell\OK$ such that $\primeq\OL$ is prime.  If $\al\in \OK$ is such that $\legendre{\al}{\primeq\OL}=1$, then $\legendre{\al}{\primeq}=1$.
\end{lemma}
\begin{proof}
Let $k_L$ be the residue field $\OL/\primeq\OL$, and $k$ be the residue field $\OK/\primeq$.  
For $z\in \OL$, let $\overline z$ denote the residue class in $k_L$. 
Since $\legendre{\al}{\primeq\OL}=1$, let $\beta\in\OL$ such that $\beta^\ell\equiv \al \mod \primeq\OL$, i.e., $\betabar^\ell=\albar$ in $k_L$. 
The degree of the extension $k(\overline \beta) / k$ is either $1$ or $\ell$ since  $k$ contains a primitive $\ell$-th root of unity. 
Since $\ell \ne p$, $k(\overline \beta)=k$. Thus, $\legendre{\al}{\primeq}=1$.\end{proof}

\begin{proposition}\label{Sec:proof:prop1}
Let $W$ be a finite subset of $\OL\backslash\set{0}$, and let $\curlyD_W$ be the set of prime numbers defined in \ref{Sec:proof:def3}.  Then $\curlyD_W$ contains a set of prime numbers in $\zz$  with positive  Dirichlet density at least
$ (p-1)/\big( \ell^{(\Sharp\, W)p!}(\ell-1) p!  \big)$.
\end{proposition}
\begin{proof}
Let $M$ be the Galois closure of $L(\sqell{\al} : \al\in W)$ over $\ratn$.
Let $M'$ be the Galois closure of $L$ over $\ratn$. Then $m:=[M':L] \not\equiv 0 \mod p$ since $f$ is defined over $\zz$, and $\deg f=:p$ is a prime number. Note that $M=M'( \sqell{\tau\al} : \al\in W\text{ and } \tau\in\GQ)$. 
Since $\zetal\in K$, and $\ell$ is a prime number, $[M: M']$ is a power of $\ell$ and, hence,  $[M:M']=\ell^n$ for some non-negative integer $n$. 
\par
It follows that $[M:K] =pm\ell^n$.    
Let $G$ denote the group $\Gal(M/\ratn)$. Then $\Gal(M/K)$ is a  subgroup of $G$, and the group $G$ contains an automorphism $\tau$ of order $p$ acting trivially on $K$.  Moreover, the subgroup $\inner{\tau}$ of $G$ is stable under conjugation since $\Gal(M/K)$ is normal in $G$. Therefore, the subset $H:=\set{\tau^k : k=1,\dots,p-1}$ is  stable under conjugation in $G$.
\par
Let $\curlyD$ be the set of all prime numbers $q$ such that $q$ is unramified in $M$ and its Frobenius automorphisms are contained in $H$.  
Then, since $[M:\ratn]$ divides $(\ell-1)\cdot p! \cdot \ell^{(\Sharp W)\, p!}$, by the Cebotarev Density Theorem, $\curlyD$ has positive Dirichlet density at least
$
       (p-1) / \big( \ell^{(\Sharp W)\, p!} (\ell-1) p! \big)
$.
Since $W$ is a finite set,  the following set of prime ideals has Dirichlet density equal to the Dirichlet density of $\curlyD$:
$$\set{ q \in \curlyD : \al \not\equiv 0 \mod \primeq \text{ for all prime ideals }\primeq \mid q\OL
        \text{ and for all } \al \in W}.$$
Thus, let us assume that $\curlyD$ is the above set.
\par
Let $q\in \curlyD$, and let $\primeQ$ be a prime ideal of  $\OM$ lying over $q$.
Let $\primeQ_L:=\primeQ\cap\OL$, and $\primeq:=\primeQ\cap \OK$.
Let us show that $\primeq\OL$ is a prime ideal. Let $f(\primeq/q)$, $f(\primeQ_L/\primeq)$,  and $f(\primeQ/\primeQ_L)$ denote the residue degrees.  Then 
\begin{equation}\label{eq:res-deg}
p=\abs{ \Frob(\primeQ/q) }=f(\primeQ/q)=f(\primeQ/\primeQ_L) f(\primeQ_L/\primeq) f(\primeq/q).
\end{equation}
Since $\tau = \Frob(\primeQ/q)\in\Gal(M/K)$, and $K/\ratn$ is Galois, $\Frob(\primeq/q)=\res_K(\tau)=1$.
Hence, $1=\abs{ \Frob(\primeq/q) } = f(\primeq/q)$.
Thus, $p=f(\primeQ/\primeQ_L) f(\primeQ_L/\primeq)$.
Since $M/L$ is Galois, $f(\primeQ/\primeQ_L)$ divides $m\ell^n\equiv 0 \mod p$ and,  hence, 
   $f(\primeQ/\primeQ_L) \not \equiv 0 \mod p$.
Therefore, $f(\primeQ/\primeQ_L)=1$ and $f(\primeQ_L/\primeq)=p$.
In other words, the prime ideal $\primeq$ remains prime in $\OL$. Moreover, $f(\primeQ/\primeQ_L)=1$ implies that $\OM / \primeQ \cong \OL / \primeQ_L$ and, hence, $\sqell{\al}$ for all $\al\in W$ are defined in $\OL / \primeQ_L$.  In other words, since $\primeQ_L=\primeq\OL$, $1=\legendre{\al }{ \primeq\OL }$ for all $\al\in W$.    
Therefore, $q\in \curlyD_W$.
\end{proof}

Recall that $\ell$ is a regular prime number.
The following lemma is the key place in our work where the additional hypothesis of $\ell$ being regular is needed.
\begin{lemma}\label{subsec:proof:lem2}
Let $W$ be a finite subset of $\OL\backslash\set{0}$ containing $\zetal$ and $\lambda:=1-\zetal$. Let $\curlyD_W$ be the set of prime numbers defined in \ref{Sec:proof:def3}. Let $\primeq$ be a place of $K$.  
If  $\primeq$ is a prime ideal of $\OK$, then we choose $\al_\primeq \in \OK$ such that $\primeq^m = \al_\primeq \OK$ where $m$ is the order of $\primeq$ in $\Cl(\OK)$, and if $\primeq=\lambda\OK$, then we choose $\alq:=\lambda$. 
\par
If $\primeq$ is an infinite place of $K$, or if $\primeq$ is a prime ideal of $\OK$ such that $\al_\primeq \in W$, then  
$\legendre{D}{\primeq}=1$
for all  positive integers $D$ supported by $\curlyD_W$.
\end{lemma}
\begin{proof} 
Suppose that  $\primeq$ is a prime ideal of $\OK$ not dividing $\ell\OK$, and $\primeq^m=\alq\OK$ where $\alq\in \OK$ and $m$ is the order of $\primeq$ in $\Cl(\OK)$. Then, since $\ell$ is regular,  $m\not\equiv 0 \mod \ell$.  Let $q$ be a prime number dividing $D$, and suppose that $\al_\primeq \in W$.  Let $q\OK$ be decomposed into $\prod_{k=1}^t \primep_j^{n}$ where $\primep_j$ are prime ideals coprime to $\ell$ and $n$ is an integer.  Since $q \in \curlyD_W$, it follows that $\legendre{\al_\primeq}{\primep_j\OL}=1$ and, hence, by Lemma \ref{lem:L-K},  $\legendre{\al_\primeq}{\primep_j}=1$. Then 
$$\legendre{\al_\primeq}{q}=\prod_{j=1}^t\legendre{\al_\primeq}{\primep_j^{n}}
                                        =\prod_{j=1}^t\legendre{\al_\primeq}{\primep_j}^{n}=1.$$
\par
Since $W$ contains $\zetal$ and $\lambda$, it follows that $\legendre{\zetal}{\primep_1} = 1$. If $\ell=2$, then $K=\ratn$, and it follows that $\legendretwo{-1}{\primep_1}=\legendretwo{2}{\primep_1}=1$.  By Lemma \ref{cor:rec1}, $\legendre{q}{\lambda\OK}=1$.
Since $\primeq \nmid \ell\OK$ and $\alq \not\equiv 0 \mod \primep_j$ for all $j=1,\dots,t$, by Corollary \ref{cor:rec2},
$$1=\legendre{\al_\primeq}{q}=\legendre{q}{\al_\primeq}=\legendre{q}{\primeq}^m.$$
Since $m\not\equiv 0 \mod \ell$, we proved that $\legendre{q}{\primeq}=1$ for all prime numbers $q$ dividing $D$.  
\par
Let $\primeq:=\lambda\OK$. Then  $\legendre{q}{\lambda\OK}=1$ for all $q \mid D$ was already shown above.
Since $D>0$, it is clear that $\legendre{D}{v}=1$ for all infinite places $v\in\MK$.
\end{proof}
\par
Let us return to the context of our superelliptic curves. 
Recall from \ref{def:context1} the point $P_1 \in X$ and the field $L:=K(P_1)$.
For all $D \in K^*$, we have the following injective maps:
\begin{equation*}
\sellam(J_D, K) \subset \HonE( K,J_D[\lambda])_{S_D} \overset{\theta^D}{\To} \thegroupell{L}.
\end{equation*}

\begin{theorem}\label{MainTheorem1-pre}
Let $K:=\ratn(\zetal)$ where $\ell$ is a regular prime number.
Let $f(x)$ be a monic polynomial of prime degree $p$ defined over $\zz$ such that $f(x)$ is  irreducible over $K$, and $\ell \ne p$.
Let $C/\ratn$  be the normalization of the \super\ $y^\ell=f(x)$.  Let $N:=\dimp\ell\sellam(J,K)$, and $M$, the number of prime ideals of $\OK$ dividing $\ell\Delta_f$.
Then there is a set $\curlyD$ of prime numbers with  Dirichlet density at least $(p-1)/\big(\ell^{(N+M+1)p!}(\ell-1) p! \big)$ such that whenever a positive integer $D$ is supported by $\curlyD$,  
$
\theta\big( \sellam(J, K) \big)=\theta^D\big( \sellam(J_D, K) \big)$.
In particular, $\dimp\ell \sellam(J, K)=\dimp\ell  \sellam(J_D, K)$.
\end{theorem}
\begin{proof}
Recall from \ref{def:C-Cq} the map $\theta : \HonE(K,J[\lambda]) \to \thegroupell{L}$.
Let $W_J$ be a subset of $\OL$ generating $\theta(\sellam(J,K))$.
For each prime ideal $\primeq$ of $\OK$ coprime to $\ell$, let us fix an element $\al_\primeq$ of $\OK$ such that 
$\al_\primeq\OK = \primeq^m$ where $m$ is the order of $\primeq$ in $\Cl(\OK)$.
For $\primeq := \lambda\OK$, we choose $\al_\primeq :=\lambda$.
\par
Recall $S_J$ from (\ref{SJ-SD}), and let
\begin{equation*}
Y_J:=\set{\zetal} \cup\ W_J
        \cup \set{\al_{\primeq}\in \OK : \primeq\in S_J \cap \MK^0}.
\end{equation*}
Let $\curlyD_{Y_J}$ be the set of prime numbers defined in \ref{Sec:proof:def3} for $W=Y_J$.
Then, by Proposition \ref{Sec:proof:prop1}, $\curlyD_{Y_J}$ contains a set $\curlyD$ of prime numbers  with Dirichlet density at least 
$(p-1)/\big( \ell^{(N+M+1)p!}(\ell-1) p! \big)$.
Since $\ell$ is regular, by Lemma \ref{subsec:proof:lem2}, if $D$ is a positive integer supported by $\curlyD$ and $\primeq \in S_J$, then
\begin{equation}\label{eq:flipit}
\legendre{D}{\primeq}=1.
\end{equation}
\par
Let $D$ be a positive integer supported by $\curlyD$.
Let us show that $\theta^D(\sellam(J_D,K)) \subset \theta(\sellam(J,K))$.
Let $\al$ be an element of $\theta^D(\sellam(J_D,K))$.  By  \ref{Sel2},
\begin{equation}\label{eq:thD-th}
\begin{aligned}
\theta(\sellam(J,K)) &= 
        \set{ \al \in \theta\big( \HonE(K,J[\lambda])_{S_J} \big) : \res_\primeq(\al) \in \Img \theta_\primeq\delta_\primeq
                \text{ for all } \primeq \in S_J };\\
\theta^D(\sellam(J_D,K)) &= 
        \set{ \al \in \theta^D\big( \HonE(K,J_D[\lambda])_{S_D} \big) : \res_\primeq(\al) \in \Img \theta^D_\primeq\delta^D_\primeq
                \text{ for all } \primeq \in S_D }.                
                \end{aligned}
\end{equation}
By Theorem \ref{Sec:proof:prop2}, we have 
\begin{equation}\label{eq:thD=th}
\theta^D\big( \HonE(K,J_D[\lambda])_{S_D} \big) = \theta\big( \HonE(K,J[\lambda])_{S_J} \big);
\end{equation}
hence, $\al$ is contained in $\theta\big( \HonE(K,J[\lambda])_{S_J} \big)$. Let $\primeq$ be a place in $S_J$, and recall the set   $S_D$ from (\ref{SJ-SD}).  Then $\primeq$ is contained in $S_D$. By Lemma \ref{hensel-lemma}, (\ref{eq:flipit}) implies that $D \in \ellthpower{\Kq}$. By Proposition \ref{prop:key2}, it follows that $\Img \theta_\primeq\delta_\primeq=\Img \theta^D_\primeq\delta^D_\primeq$ and, hence, $\al$ is contained in $\Img \theta_\primeq\delta_\primeq$ since $\al \in \theta^D(\sellam(J_D,K))$ and $\primeq \in S_D$.  Thus, $\al \in \theta(\sellam(J,K))$.
\par
Let us show that $\theta(\sellam(J,K)) \subset \theta^D(\sellam(J_D,K)) $. Let $\al$ be an element of $\theta(\sellam(J,K))$.
By (\ref{eq:thD=th}), $\al$ is contained in $\theta^D\big( \HonE(K,J_D[\lambda])_{S_D} \big)$. If $\primeq$ is a place in $S_J$, then by (\ref{eq:flipit}), $D \in \ellthpower{\Kq}$, and by Proposition \ref{prop:key2}, $\Img \theta_\primeq\delta_\primeq=\Img \theta^D_\primeq\delta^D_\primeq$.  Thus, $\al \in \Img \theta^D_\primeq\delta^D_\primeq$ for all $\primeq \in S_J$. Let $\primeq$ be a prime ideal in $S_D\minuS S_J$.
Since $D$ is supported by $\curlyD$, the prime ideal $\primeq$ is contained in $\curlyD$ and, hence, $\legendre{\beta}{\primeq}=1$ for all $\beta \in W_J$.  Since $W_J$ generates $\theta(\sellam(J,K))$, it follows that $\legendre{\al}{\primeq\OL}=1$, i.e., 
$\al \in \ellthpower{L_\primeP}$ where $\primeP:=\primeq\OL$. Recall that $\res_\primeq$ is a map: 
$\thegroupell{L} \to \Cq:=\prod_{P\in Z(J,\primeq)} \thegroupell{\Kq(P)}$. Since $\primeq\OL$ is prime, $f(x)$ is irreducible over $\Kq$ and, hence, $Z(J,\primeq)$ contains a single point $P$. Thus, $\Cq=\thegroupell{ L_\primeP}$.  Thus, $\res_\primeq(\al)=1$ and, in particular, $\res_\primeq(\al) \in \Img \theta^D_\primeq\delta^D_\primeq$. We established that $\res_\primeq(\al) \in \Img \theta^D_\primeq\delta^D_\primeq$ for all $\primeq \in S_D$.  Therefore, $\al$ is contained in 
$\theta^D(\sellam(J_D,K))$. Since $\theta$ and $\theta^D$ are injective maps, the dimensions of the two selmer groups are equal to each other.\end{proof}

Recall that $\curlyP_\ell(X)$ is the set of positive $\ell$-th power free integers up to $X$.
The following lemma is easily deduced from \cite{serre:1976}, Theorem 2.4, and the proof is left to the reader:
\begin{lemma}\label{prop:multi-set}
Let $\curlyD$ be a set of prime numbers in $\zz$ with positive Dirichlet density $<1$.
Let $\curlyP_\ell(\curlyD,X)$ denote the set of all positive integers up to $X$, which are $\ell$-th power-free and supported by $\curlyD$.
Let $\curlyP_\ell(\curlyD):=\cup_{X=1}^\infty \curlyP_\ell(\curlyD,X)$. 
Then there are positive constants $c$ and $\ep<1$ depending on $\curlyP_\ell(\curlyD)$ such that 
              $\SharP \curlyP_\ell(\curlyD,X) \sim  \frac{ X }{ (\log X)^\ep }$ as $X \to \infty$.
\par
Let $D_0$ be a positive $\ell$-th power free integer.  Then  there is a positive constant  $\ep<1$ depending on $\curlyD$ such that 
\begin{equation*}
        \SharP \set{D_0\,D \in \curlyP_\ell(X) : D\in\curlyP_\ell(\curlyD,X) } \gg \frac{ X }{ (\log X)^\ep }.
        \end{equation*}
\end{lemma}

\begin{theorem}\label{MainTheorem1}
Let $K$, $f(x)$, and $C/\ratn$ be as in Theorem \ref{MainTheorem1-pre}.
Let $D_0$ be a positive integer. Let $N:=\dimp\ell\sellam(J_{D_0},K)$, and let $M$ be the number of prime ideals of $\OK$ dividing $\ell\Delta_f\, D_0$.
Then there is a set $\curlyD$ of prime numbers with  Dirichlet density at least $(p-1)/\big( \ell^{(N+M+1)p!}(\ell-1) p! \big)$ such that whenever a positive integer $D$ is supported by $\curlyD$, 
\begin{equation}
\dimp\ell\Sel{\lambda}(J_{D_0\,D}, K) = N.\label{eq:MT1'}
\end{equation}
Moreover, there is a positive constant $\ep<1$ such that 
\begin{equation}
\Sharp\ \set{ \DPlX :\sellam(J_{D},K)=N} \gg_{J,D_0} \frac{ X }{ (\log X)^{\ep} }.\label{eq:MT1}
\end{equation}
\end{theorem}
\begin{proof}
Recall  that $J_{D_0}$ is the Jacobian variety of the normalization of $y^\ell = (D_0)^p\ f(x/D_0)$.  Let $g:=(D_0)^p\ f(x/D_0)$. Then $\Delta_g=\Delta_f\, D_0^{b}$ for some positive integer $b$ and, hence, the number of prime ideals of $\OK$ dividing $\ell \Delta_g$ is equal to $M$.
Since $(J_{D_0})_D=J_{D_0\,D}$, by applying  Theorem \ref{MainTheorem1-pre} to $J_{D_0}$, (\ref{eq:MT1'}) is proved. The proof of (\ref{eq:MT1}) follows immediately from Lemma \ref{prop:multi-set}. \end{proof}
\begin{corollary}\label{cor:rank<=n}
Assume the same hypotheses in Theorem \ref{MainTheorem1}.  
Then there is a positive constant $\ep<1$ such that 
\begin{equation}
\Sharp\ \set{ \DPlX :\rank J_D(\ratn)\le N} \gg_{J,D_0} \frac{ X }{ (\log X)^{\ep} }.\label{eq:MT2}
\end{equation}
\end{corollary}
\begin{proof}
Note $\sellam(J_D,K) \cong J_D(K)/ \lambda J_D(K) \oplus \sha(J_D,K)[\lambda]$. Then, by \cite{schaefer:1998}, Corollary 3.7 and Proposition 3.8, $\sellam(J_D,K) \ge \rank J_D(\ratn)$. Hence, (\ref{eq:MT2}) follows. \end{proof}

\begin{corollary}\label{cor:stoll}
Let $C/\ratn$ be the normalization of a hyperelliptic curve $y^2=f(x)$ where $f(x)\in\zz[x]$ is monic and irreducible over $\ratn$ and has odd prime degree $p\ge 5$.
Suppose that there is a positive integer $D_0$ such that $N:=\SEL(J_{D_0},\ratn) < (p-1)/2$. 
Then there is a positive constant $\ep<1$ depending on $C$ and $D_0$  such that 
\begin{equation*}
\Sharp\ \set{ \DPtwoX : \Sharp\ C_D(\ratn) \le 2N+1 }
                \gg \frac{X}{(\log X)^\ep}.
\end{equation*}
\end{corollary}
\begin{proof}
Let $\iota$ be \textit{the hyperelliptic involution on $C_D$.} 
In \cite{stoll:2004}, Theorem 1.1, Stoll proved that if $D$ is coprime to a fixed finite set $T$ of  prime numbers determined by $C$ and $D_0$, then 
any set $S\subset C_D(\ratn)$ such that $\Sharp\ S \le (p-1)/2$ and $S\cap \iota(S)=\emptyset$ generates
a subgroup of rank $\Sharp S$ in $J_D(\ratn)$.
\par
By Theorem \ref{MainTheorem1}, 
\begin{equation*}
\Sharp\ \set{ \DPtwoX : \dimp2\Sel{2}(J_D,\ratn)=N }
                \gg \frac{X}{(\log X)^\ep}.\label{eq:stoll-count}
\end{equation*}
We count the $D$'s supported by a set $\curlyD$ of prime numbers with positive Dirichlet density in order to find such a lower bound. Since $T$ is finite, we may assume that $\curlyD$ does not intersect $T$.
Let $D $ be a positive integer not supported by $T$ such that $\dimp2\Sel{2}(J_D ,\ratn)=N$. Then, by Stoll's theorem, 
Thus, $C_D (\ratn)$ can not contain a subset $S$ such that $\Sharp\, S > N$ and $S \cap \iota(S) = \emptyset$; otherwise, $\rank J_D(\ratn) > N$.
Since $C_D (\ratn)$ does not contain a point fixed under the involution, except $\infty$, we conclude 
$\Sharp\,C_D (\ratn) \le 2N+1$. \end{proof}
\begin{corollary}\label{cor:cubic}
Let $E/\ratn$ be an elliptic curve given by $y^2=x^3-A$ where $A$ is a positive square-free integer such that 
$A\equiv 1$ or $25 \mod 36$ and $\dimp3 \Cl(\ratn(\sqrt{-A}))[3]=0$.  For a non-zero cube-free integer $D$, let $E_D$ be the cubic twist: $y^2=x^3-A\,D^2$.  Then there is a positive integer $\ep <1$ such that 
\begin{equation*}
 \Sharp\, \set{ D \in \curlyP_3(X) : \rank E_D(\ratn) = 0 } \gg \frac{ X }{ (\log X)^\ep }.
 \end{equation*}
 \end{corollary}
 \begin{proof}
 By \cite{stoll:1998}, Corollary 2.1, $\dimp3 \sellam(E,K)=0$ where $\lambda=1-\zeta_3$ and $K=\ratn(\zeta_3)$. 
 As $A$ is square-free and coprime to $3$, the polynomial $y^2+AD^2$ is irreducible over $K$. 
 The result
 follows immediately from Theorem \ref{MainTheorem1} with $\ell=3$. 
 \end{proof}

\subsection{The case of function fields}
In this section, we shall prove two results for function fields that are analogous to Theorem \ref{MainTheorem1}.
For Lemma \ref{lem:LtoK}, Proposition \ref{prop:pos-Dirch-FF}, and Lemma \ref{prop:Dg=1}, let us assume a slightly more general context:
Let $k$ be a finite field containing a primitive $\ell$-th root of unity $\zetal$; hence, $\Char k \ne \ell$.
Let $K/k$ be a \globalfieldv\ such that $\Cl(\OK)[\ell]\not\equiv 0 \mod \ell$ defined in \ref{rel-infty}. Let $L$ be a field extension of $K$ such that $p:=[L:K]$ is a prime number not equal to $\ell$.

\begin{definition}\label{DW-FF}
Let $W$ be a finite subset of $\OL\minuS\set{0}$.  We denote by $\curlyD_W$ the set of  prime ideals $\primeq$ of $\OK$ which satisfy the following properties: 
\begin{enumerate}
\item The ideal $\primeq\OL$ is prime;
\item for all $\al \in W$, $\al \not\equiv 0 \mod \primeq\OL$, and $\legendre{\al}{\primeq\OL} =1$.
\end{enumerate}
\end{definition}

\begin{lemma}\label{lem:LtoK}
Let $\primeq$ be a prime ideal of $\OK$ such that $\primeq\OL$ is prime. 
If $\al \in \OK$ such that $\legendre{\al}{\primeq\OL}=1$, then $\legendre{\al}{\primeq}=1$.
\end{lemma}
\begin{proof}
The proof is similar to that of Lemma \ref{lem:L-K}.\end{proof}

\begin{proposition}\label{prop:pos-Dirch-FF}
Let $W$ be a finite subset of $\OL\backslash\set{0}$.  Then $\curlyD_W$ contains a set of prime ideals of $\OK$ with positive  Dirichlet density at least $(p-1)/\big( \ell^{(\Sharp W)p!} p! \big)$.
\end{proposition}
\begin{proof}
The proof is similar to that of Proposition \ref{Sec:proof:prop1}.
\end{proof}

\par
Recall that $K/k$ is a \globalfieldv, and let $\Oinf$, $\Minf$, and $\piinf$ be the discrete valuation ring, the maximal ideal, and a uniformizer at $v_\infty$, respectively. 

\begin{definition}\label{at-infty}
Let $\Kinf$ be the completion of $K$ at $v_\infty$. Let us define
$$ \legendre{g}{v_\infty}:=\begin{cases}
                1 &\text{ if } g \in (\Kinf^*)^\ell\\
                -1&\text{ if } g \not\in (\Kinf^*)^\ell
                \end{cases}.
                $$
Note that each non-constant element $g$ of $\OK$ is not an element of the valuation ring $\OK[\infty]$, i.e., $\ord_{v_\infty}(g)<0$. 
\textit{The leading coefficient of a nonzero element $g$ of $\OK$} is the constant $a \in k^*$ such that $\piinf^m\, g \equiv a \mod \piinf$ for some $m \in \zz$. 
The element $g$ is \textit{monic} if the leading coefficient is $1$.
\textit{The degree of an element $g$ of $\OK$}, denoted by $\deg(g)$, is $-\ordinf(g)$.
Let $g$ be an element of $\OK$ with the leading coefficient $a \in k^*$. Then $g \in \big( \Kinf^* \big)^\ell$ if and only if 
$a \in (k^*)^\ell$ and $\deg(g)$ is divisible by $\ell$.
\end{definition}

\begin{lemma}\label{prop:Dg=1} 
Let $W$ be a finite subset of $\OL\backslash\set{0}$ containing $-1 \in k$, and let $\curlyD_W$ be the set of prime ideals of $\OK$ defined in \ref{DW-FF}.  Let $\primeq$ be a prime ideal of $\OK$, and $\alq$, an element of $\OK$ such that 
$\alq\OK=\primeq^m$ where $m$ is the order of $\primeq$ in $\Cl(\OK)$.
If  $\alq \in W$, and $D$ is a nonzero monic element in $\OK$, supported by $\curlyD_W$ such that $\deg(D) \equiv 0 \mod \ell$, then 
$$ \legendre{D}{v_\infty}=\legendre{D}{\primeq}=1.$$
\end{lemma}
\begin{proof} 
Suppose that $\alq \in W$, and let $D$ be a nonzero monic element in $\OK$, supported by $\curlyD_W$ such that $\deg(D) \equiv 0 \mod \ell$.
Then, $\legendre{D}{v_\infty}=1$. 
Following the proof of Lemma \ref{subsec:proof:lem2}, we find $\legendre{\alq}{D}=1$.
Since $-1 \in W$, it follows that $\legendre{-1}{D}=1$.
By Lemma \ref{lem:rec2}, $1=\legendre{\alq}{D}=\legendre{D}{\alq}$.  It follows that 
$1= \legendre{D}{\alq}=\legendre{D}{\primeq}^m$ where  $m \not\equiv 0 \mod \ell$.
Therefore, $\legendre{D}{\primeq}=1$.\end{proof}

\begin{theorem}\label{prop:hyper-FF-no2}
Assume that $\Sharp\ \Cl(\OK) \not\equiv 0 \mod \ell$.
Let $f(x)$ be a monic polynomial of prime degree $p$ defined over $\OK$ such that $f(x)$ is  irreducible over $K$, and $\ell \ne p$.
Let $C/K$  be the normalization of the \super\ $y^\ell=f(x)$.  Let $N:=\dimp\ell\sellam(J,K)$, and let $M$ be the number of prime ideals of $\OK$ dividing $\ell\Delta_f$.
Then there is a set $\curlyD$ of prime ideals with  Dirichlet density at least $(p-1)/(\ell^{(N+M+1)p!}p! )$ such that whenever $D$ is a monic element of $\OK$ supported by $\curlyD$ such that $\deg(D)$ is divisible by $\ell$,
\begin{equation*}
\theta\big( \sellam(J, K) \big)=\theta^D\big( \sellam(J_D, K) \big).
\end{equation*}
\end{theorem}
\begin{proof} 
Recall $S_J$ from \ref{SJ-SD}.
For each prime ideal $\primeq$ of $\OK$, choose a monic element $\al_\primeq$ of $\OK$ such that $\alq \OK = \primeq^m$ where $m$ is the order of $\primeq$ in $\Cl(\OK)$.
Let $W_J$ be the subset of $\OL$ generating $\theta(\sellam(J,K))$. Let 
$$Y_J:=\set{-1}\cup W_J \cup\set{ \alq \in \OK : \primeq \in S_J \cap \MK^0},$$
and let $\DYJ$ be the set of prime ideals defined in \ref{DW-FF} for $W=Y_J$.
Then, by Proposition \ref{prop:pos-Dirch-FF}, $\DYJ$ contains a set of prime ideals $\curlyD$ with  Dirichlet density at least $(p-1)/(\ell^{(N+M+1)p!} p!)$.
The proof of $\theta\big( \sellam(J, K) \big)=\theta^D\big( \sellam(J_D, K) \big)$ is identical with that of Theorem \ref{MainTheorem1}, when Lemma \ref{prop:Dg=1} is applied to $W=Y_J$. 
\end{proof}

\begin{corollary}\label{cor:hyper-FF-no2}
Assume the same hypotheses as in Theorem \ref{prop:hyper-FF-no2}.
Let $D_0$ be a nonzero element of $\OK$. Let $N:=\dimp\ell\sellam(J_{D_0},K)$, and $M$, the number of prime ideals of $\OK$ dividing $\ell\Delta_f\, D_0$ . Then there is a set $\curlyD$ of prime ideals of $\OK$ with  Dirichlet density at least $(p-1)/(\ell^{(N+M+1)p!} p! )$ such that whenever $D$ is a monic element of $\OK$ supported by $\curlyD$ such that $\deg(D)$ is divisible by $\ell$,
$$\dimp\ell\sellam( J_{D_0\,D},K)=N.$$
\end{corollary}
\begin{proof} 
The proof is similar to that of Theorem \ref{MainTheorem1}.\end{proof}

\begin{definition}\label{def:manyD}
Given a set $\curlyD$ of prime ideals of $\OK$, there are indeed infinitely many classes in $\thegroupell{K}$ represented by $D \in \OK$ supported by $\curlyD$ and $\deg(D)\equiv 0\mod \ell$.
Let $\set{\primep_1,\cdots,\primep_n}$ be a finite set of prime ideals of $\OK$ contained in $\curlyD$ for $n\ge 2$.
Let $\primep_i^{m_i}=\al_{\primep_i}\OK$ for some $\al_{\primep_i} \in \OK$ where $m_i$ is the order of $\primep_i$ in $\Cl(\OK)$.
Then, since $\Cl(\OL) \not\equiv 0 \mod \ell$ and $\deg(\al_{\primep_n})=\ord_{\primep_n}(\al_{\primep_n})$, there is a positive integer $s$ such that $D:= \al_{\primep_n}^s \prod_{i=1}^{n-1} \alpi$ has degree divisible by $\ell$.  
\end{definition}

\begin{theorem}\label{thm:functionfield1}
Assume the same hypotheses as in Theorem \ref{prop:hyper-FF-no2}.
Suppose that $f(x)$ is defined over $k$. Let $k'$ be the finite extension of $k$ of degree $p:=\deg(f)$. Let $L:=K\otimes k'$, and $\OL$, the integral closure of $K$ in $L$.
Suppose that $\dimp\ell \Cl(\OL)[\ell]=0$. Then $\dimp\ell\sellam(J,K)=0$. 
\par
Let $E/k$ be the \textit{constant} Jacobian variety of the normalization of the superelliptic curve $y^\ell = f(x)$ over $k$.
Then there is a set $\curlyD$ of prime ideals of $\OK$ with Dirichlet density $(p-1) / p$ such that 
whenever $D$ is an element of $\OK$ supported by $\curlyD$,
$$\dimp\ell\sellam(J_D,K)=0\quad \text{ and }\quad
\Sharp\ C_D(K) \le \Sharp\ E(k).
$$
\end{theorem}
\begin{proof}
Note that $L/K$ is Galois.
Recall $S_J$ from (\ref{SJ-SD}). Then $S_J=\MKinf$, and
$\HonE(K,J[\lambda])_{\MKinf}\cong \ker\big( \normLK : L(\MKinf,\ell) \to K(\MKinf,\ell) \big)$.
Note that 
$\Sharp\,\sellam(J,K) \le \Sharp\,  \HonE(K,J[\lambda])_{\MKinf}
                        =\Sharp\, \ker\big( \normLK : L(\MKinf,\ell) \to K(\MKinf,\ell) \big)$.
Since the subgroup $\Cl(\OL)[\ell]$ is trivial,        
by Lemma \ref{lem:ker-triv} below, $\HonE(K,J[\lambda])_{\MKinf}=1$ and, hence, $\sellam(J,K)=0$.
\par
Let $\curlyD$ be the set of prime ideals $\primeq$ of $\OK$ such that $\primeq\OL$ is prime.  Since $\Gal(L/K)$ has order $p$, by the Cebotarev Density Theorem, $\curlyD$ has Dirichlet density $(p-1)/p$.  Let $D$ be an element of $\OK$ supported by $\curlyD$.
Then, by Theorem \ref{Sec:proof:prop2}, $\theta^D(\HonE(K,J_D[\lambda])_{S_D})=\theta(\HonE(K,J[\lambda])_{S_J})=1$ and, hence,  $\dimp\ell\sellam(J_D,K)=0$. 
\par
By \cite{schaefer:1998}, Corollary 3.7,
$
\rank J_D(K) \le (\ell-1)\dimp\ell\sellam(J_D,K)=0$.
Suppose that $D$ is not a unit in $\OK$, and let $F:=K(\sqrt[\ell]{D})$. Then $k$ is algebraically closed in $F$, and $F \subset K\sep$ since $\ell \ne \Char K$.
Since $(C_D)_F \cong C_F$ as $F$-schemes, by Proposition \ref{prop:k-to-K} below, $\Sharp\ J_D(K)\tor \le \Sharp\ E(k)$.
Since $J_D(K) = J_D(K)\tor$, and $C_D(K) \injects J_D(K)$, it follows that 
$ \Sharp\ C_D(K) \le \Sharp\ J_D(K)\tor \le \Sharp\ E(k)$.
\end{proof}
\par

\begin{lemma}\label{lem:ker-triv}
(Suppose that $\MKinf$ consists of a single place of degree $1$.)  Let $L := K\otimes k'$ where $k'$ is a finite (separable) extension of $k$ of degree $d$ coprime to $\ell$, and suppose that $\Cl(\OL)[\ell]$ is trivial.
Then 
$ \ker \big( \normLK : L(\MKinf,\ell) \to K(\MKinf,\ell) \big)=1$.
\end{lemma}
\begin{proof}
Since the extension $L/K$ is obtained from base change, $M_L^\infty$ still consists of a single place of degree $1$,
and $L$ is a function of field of one variable over $k'$.
With $S=\MK^\infty$, $L(\MKinf,\ell) \cong \thegroupell{\OL}$ since $\Cl(\OL)[\ell]=1$. 
It follows from Proposition \ref{prop:dino} below that $\OL^*=(k')^*$. Hence, 
\begin{equation}\label{eq:kernormfinite}
\ker\big( \normLK : L(\MKinf, \ell) \To K(\MKinf,\ell) \big) 
        \cong \ker\big( \Norm_{k'/k} : \thegroupell{k'} \To \thegroupell{k} \big).
\end{equation}        
Let $a$ be a positive integer such that $da\equiv 1 \mod \ell$. Note that there is a section $u : \thegroupell{k} \to \thegroupell{k'}$ to the norm map $\Norm_{k'/k}$ in (\ref{eq:kernormfinite}), given by $\cls(\al) \mapsto \cls(\al^a)$. 
Hence, the norm map $\Norm_{k'/k}$ is surjective. 
Since $k^*$ and $(k')^*$ are both cyclic groups of order divisible by $\ell$, it follows that $\Norm_{k'/k}$ is an isomorphism of the $\finiteell$-vector spaces. Therefore, the kernel is trivial.\end{proof}

\begin{proposition}\label{prop:dino}
Let $k'$ be a finite field.
Let $\curlyZ/k'$ be a smooth complete curve with function field $F$ such that $\curlyZ$ has a rational divisor $v_\infty$.
Then   there is a $k'$-morphism: $\curlyZ \to \Pone_{k'}$ such that $\vinf$ is totally ramified over a rational divisor in $\Pone_{k'}$. 
\par
Let $\OF$ be the ring of integers defined in \ref{rel-infty} with choice of $M_F^\infty:=\set{v_\infty}$.
Then the group of units $\OF^*$ is  $(k')^*$, and 
$
\Sharp\, \Cl(\OF) = \Sharp\, \Pico(\curlyZ)
                        $.
Hence, the class number of $\OF$ does not depend on the choice of a rational divisor on $F$.
\end{proposition}
\begin{proof} 
Using the Riemann-Roch theorem, we can find a function $g$ in $k'(\curlyZ)$ with poles supported only by $v_\infty$. Then the function $g$ induces a morphism $\curlyZ \to \Pone_{k'}$ such that $\vinf$ is totally ramified over a rational point $\infty \in \Pone_{k'}$.
To finish the proof, use \cite{dino:1996}, Sec VIII, p.~299-300.\end{proof}

\begin{proposition}\label{prop:k-to-K}
Let $k$ be a perfect field, and let $K$ be a field extension of $k$ such that $k$ is algebraically closed in $K$.
Let $E/k$ be a smooth complete geometrically connected curve with a $k$-rational point. 
\par
Let $C'/K$ be a twist of $E_K/K$, and suppose that there is a field extension $F$ of $K$ such that $k$ is algebraically closed in $F$ and such 
that $C'_F \cong E_F$ (as $F$-schemes). Let $J/k$ be the Jacobian variety of $E/k$, and let $J_{C'}/K$ be the Jacobian variety of $C'/K$.
Then 
        $ J_{C'}(K)\tor \injects J(k)
        $.
\end{proposition}
\begin{proof} 
Note that 
$J(\kbar)\tor = J(\Kbar)\tor = J(K\sep)\tor$. Then,
$J_{C'}(K)\tor \subset J_{C'}(F)\tor \cong J_F(F)\tor \cong J(F)\tor = J(k)$.
 \end{proof}

\section{The Jacobian Varieties With $\lambda$-Torsion Points}\label{sec:full-torsion}

In this section we prove Theorem \ref{MainTheorem3} and its analogue for a function field, introduced in Section \ref{sec:intro}.
Let $\ell$ be a prime number, and let $K:=\ratn(\zetal)$ for which we assume $\ell$ is regular, or a \globalfieldv\ such that $\Sharp\,\Cl(\OK) \not \equiv 0 \mod \ell$.
Let $f(x)$, $C/K$, $C_D/K$, $J/K$, and $J_D/K$ be as in Section \ref{sec:def-prop}, and we keep the notation used in that section.
Let $\Delta_f$ be the discriminant of $f$.
Recall  $z_1,\dots,z_d \in K\sep$, the roots of $f(x)$, and suppose that $z_d$ is contained in $ K$. 
 Recall that $Z(J):=\set{T_1,\dots,T_s}$ is a  set of representatives of $\GGK$-orbits in $X$, and  $L_i:=K(T_i)$ for $i=1,\dots,s$. Let $L$ be the compositum of $L_1,\dots,L_s$ in $K\sep$.

\begin{definition}\label{def:DYJp}
Let us fix a set of generators of $\theta(\sellam(J,K))$, and note that each generator is an $s$-tuple with entries in $L$. 
Let $W_J$ be the union of all entries of the generators. Then $W_J$ is a subset of $L^*$. 
For each prime ideal $\primeq$ of $\OK$, choose an element $\al_\primeq$ of $\OK$ as in the proof of Theorem \ref{MainTheorem1-pre} or  Theorem \ref{prop:hyper-FF-no2}, depending on the cases of $K$.
Recall $S_J:=\MK^\infty\cup\set{ \primeq \in \MKfinite : \primeq \mid \ell\Delta_f\OK }$, and let
\begin{equation*}
Y_J:=\set{\zetal,-1} \cup\ W_J
        \cup \set{\al_{\primeq}\in \OK : \primeq\in S_J \cap \MK^0}.
\end{equation*}
When $K=\ratn(\zetal)$, let $M$ be the Galois closure of $L(\sqrt[\ell]{\al} : \al\in Y_J)$ over $\ratn$. When $K$ is a function field, 
let $M$ be the Galois closure of $L(\sqrt[\ell]{\al} : \al\in Y_J)$ over $K$.
\par
If $K=\ratn(\zetal)$, let us denote by $\DYJJ$ the set of prime numbers $q$ in $\zz$ such that $q$ splits completely in $\OM$ and coprime to $\al$ for all $\al \in Y_J$. If $K$ is a function field, let us denote by $\DYJJ$ the set of prime ideals $\primeq$ of $\OK$ such that $\primeq$ splits completely in $\OM$ and  coprime to $\al$ for all $\al \in Y_J$.
By the Cebotarev Density Theorem with $H$ being the trivial subgroup of $\Gal(M/\ratn)$ or $\Gal(M/K)$, depending on both cases of $K$, the set of prime ideals of $\zz$ or $\OK$ that split  completely in $M$ has positive Dirichlet density.  Since there are finitely many prime ideals $\primeq$ of $\zz$ or $\OK$ that are not coprime to $\al$ for some $\al \in Y_J$, $\DYJJ$ has positive Dirichlet density.
\end{definition}

\par
Recall from \ref{at-infty} the definition of $\legendre{\bullet}{v_\infty}$ when $K$ is a function field, and that
$\Sharp\ \Cl(\OK) \not\equiv 0 \mod \ell$.
\begin{lemma}\label{prop:DYJJ}
Suppose that $K=\ratn(\zetal)$.
If $D$ is  a positive integer supported by $\DYJJ$, and $\primeq$ is a place in $S_J$, then 
$\legendre{D}{\primeq}=1$.  
\par
Suppose that $K$ is a function field.
Let $D$ be  a monic element of $\OK$, of degree  divisible by $\ell$  supported by $\DYJJ$.  If $\primeq$ is a place in $S_J$, then 
$\legendre{D}{\primeq}=1$.  
\end{lemma}
\begin{proof}
The proof is identical with those of Lemma \ref{subsec:proof:lem2} and \ref{prop:Dg=1}.\end{proof}

\begin{proposition}\label{prop:Sel-infinity}
Suppose that $K=\ratn(\zetal)$.  Let $D$ be a positive $\ell$-th power free integer in $\zz$ which is supported by $\curlyD'_{Y_J}$.
Then 
\begin{equation}\label{eq:ineq1}
 \dimpell\sellam(J_D,K)> \dimpell\sellam(J,K).
 \end{equation}
\par
Suppose that $K$ is a function field.
Let $D$ be a monic element of $\OK\minuS\set{0}$ of degree divisible by $\ell$ such that $D\OK$ \notellthpower\ and $D$ is supported by $\curlyD'_{Y_J}$.  Then 
\begin{equation}\label{eq:ineq2}
 \dimpell\sellam(J_D,K)> \dimpell\sellam(J,K).
 \end{equation}
For both cases of $K$, 
$\lim\sup_D\ \dimp\ell\Sel{\lambda}(J_D,K) =\infty$.
\end{proposition}
\begin{proof}
Let $D$ be  a positive $\ell$-th power free integer in $\zz$ which is supported by $\curlyD'_{Y_J}$ or 
a monic element $D$ of $\OK$ of degree divisible by $\ell$ supported by $\DYJJ$.
Let $S_D:=S_J\cup\set{\primeq \in \MKfinite : \primeq \mid D\OK}$.  Then,
\begin{equation*}
\theta( \HonE(K,J[\lambda])_{S_J} ) = \prod_{i=1}^s L_i(S_J,\ell),\quad
\theta^D( \HonE(K,J_D[\lambda])_{S_D} ) = \prod_{i=1}^s L_i(S_D,\ell).
\end{equation*}
By \ref{Sel2}, 
\begin{equation}\label{eq:Sel-descr}
\begin{aligned}
\theta(\sellam(J,K)) &= 
        \set{ \al \in \theta\big( \HonE(K,J[\lambda])_{S_J} \big) : \res_\primeq(\al) \in \Img \theta_\primeq\delta_\primeq
                \text{ for all } \primeq \in S_J };\\
\theta^D(\sellam(J_D,K)) &= 
        \set{ \al \in \theta^D\big( \HonE(K,J_D[\lambda])_{S_D} \big) : \res_\primeq(\al) \in \Img \theta^D_\primeq\delta^D_\primeq
                \text{ for all } \primeq \in S_D }.
\end{aligned}
\end{equation}
\par
Let us show that $\theta(\sellam(J,K)) \subset \theta^D(\sellam(J_D,K))$. Let $\al$ be an element of $\theta(\sellam(J,K))$. Then $\al \in \prod_{i=1}^s L_i(S_J,\ell)$. Since $S_J$ is contained in $S_D$, we have $\al \in \theta^D\big( \HonE(K,J_D[\lambda])_{S_D} \big) = \prod_{i=1}^s L_i(S_D,\ell)$.
If $\primeq$ is a place in $S_J$, then by Lemma \ref{prop:DYJJ}, $D$ is an $\ell$-th power in $\Kq$. By Proposition \ref{prop:key2}, 
$\Img \theta^D_\primeq\delta^D_\primeq=\Img \theta_\primeq\delta_\primeq$ and, hence, $\al \in \Img \theta^D_\primeq\delta^D_\primeq$. Let $\primeq$ be a prime ideal in $S_D \minuS S_J$.
Then $\primeq \mid D$ and, hence, $\primeq \in \curlyD'_{Y_J}$.  Suppose that $\al$ is represented by $(\beta_1,\dots,\beta_s)$ for some $\beta_i \in \OL$. Then, since $Y_J$ contains $W_J$ defined in \ref{def:DYJp}, and all $\sqrt[\ell]{\beta_i}$ are defined over $\Kq$, it follows  $\res_\primeq(\al)=1 \in \Cq$ and, in particular, $\res_\primeq(\al) \in \Img \theta^D_\primeq\delta^D_\primeq$.
Therefore, by (\ref{eq:Sel-descr}), $\al$ is contained in $\theta^D(\sellam(J_D,K))$.
\par
Note that $\theta(\sellam(J,K)) \subset \prod_{i=1}^s L_i(S_J,\ell)$. Lemma \ref{lem:2-torsion} below shows that there is an element $\al$ of $\theta^D(\sellam(J_D,K))$ which is  contained in $\prod_{i=1}^s L_i(S_D,\ell) \minuS \prod_{i=1}^s L_i(S_J,\ell)$.  Therefore, $\al \not\in \theta(\sellam(J,K))$, and we proved that $\dimp\ell\sellam(J_D,K)>\dimp\ell\sellam(J,K)$.
\par 
If $K$ is a function field, then as illustrated in \ref{def:manyD},  there is a monic element $D$ of $\OK$  of degree divisible by $\ell$ such that $D\OK$ \notellthpower\ and $D$ is supported by $\DYJJ$.
 Using induction, one can easily show that
  $\lim\sup_D\ \dimp\ell\Sel{\lambda}(J_D,K) =\infty$.
\end{proof}
\par
Recall  that $P_d^D:=[(z_d\, D,0)-(\infty)]$ is a point in $J_D[\lambda](K\sep)$ but $P_d^D \not \in \XD$.
\begin{lemma}\label{lem:2-torsion}
Let $D$ be an $\ell$-th power free positive integer supported by $\DYJJ$, or a monic element of $\OK$ supported by $\DYJJ$ such that $D\OK$ \notellthpower.
Consider the map
$$\theta^D\delta^D : \WMDell \To \curlyC:=\prod_{i=1}^s \thegroupell{L_i}.$$
Then the $K$-rational point $P_d^D\in J_D[\lambda](K)$ is mapped to 
$\pprod L_i(S_D,\ell) \minuS \pprod L_i(S_J,\ell)$ under $\theta^D\delta^D$.
\end{lemma}
\begin{proof}
Since $D$ is supported by $\DYJJ$, $D$ is coprime to $\Delta_f$ and to all prime ideals $\primeq\in S_J$.
Recall that $P_i^D:= [(z_i\, D,0)-(\infty)] $ for $i=1,\dots,d$.
Since $P_d^D:=[(z_d\, D,0)-(\infty)]$ is a point in $J_D(K)$,  by Lemma \ref{cor:coboundary},
\begin{multline}
\theta^D(\delta^D( P_d^D )) = ( f_{P^D}(z_d\,D,0) : P \in Z(J))
        =( z_d\,D - D\ x(P) : P \in Z(J))\\
        =( D(x(P_d) - x(P) ) : P \in Z(J)).
        \end{multline}
Note that for all $P\in Z(J)$, the difference $x(P_d)-x(P)$ divides $\Delta_f$ and hence, it is coprime to $D$.
  Since 
$D\OK$ \notellthpower, it follows that 
$D\,\big( x(P_d) - x(T_i) \big) \in L_i(S_D,\ell) \minuS L_i(S_J,\ell)$
for all $i=1,\dots,s$.
\end{proof}
 \begin{theorem}\label{MainTheorem3}
Suppose that $K=\ratn(\zetal)$.
Given a positive integer $n$, there is a positive constant $\ep<1$ depending on $C$ and $n$ such that 
\begin{equation}\label{eq:main2}
\Sharp\,\set{ \DPlX : \dimp\ell\Sel{\lambda}(J_D, K) > n }
                \gg \frac{X}{(\log X)^\ep}.
\end{equation}
\end{theorem}
\begin{proof}
By Proposition \ref{prop:Sel-infinity}, there is a positive $\ell$-th power-free rational integer $D_0$ such that  \\   
$\dimpell\sellam( J_{D_0},K) > n$.
Proposition \ref{prop:Sel-infinity} applied to $J=J_{D_0}$ and Lemma \ref{prop:multi-set} together imply (\ref{eq:main2}).
\end{proof}

\begin{theorem}\label{MainTheorem3'}
Suppose that $K$ is a function field.
Given a positive integer $n$, there are $D_0 \in \OK$ and  a set of prime ideals $\curlyD$ of $\OK$ with positive Dirichlet density such that whenever $D$ is
a monic element  of $\OK$ of degree divisible by $\ell$ such that $D\OK$ \notellthpower\ and $D$ is supported by $\curlyD$,
\begin{equation}\label{eq:main2'}
\dimp\ell\Sel{\lambda}(J_{D_0\,D}, K) > n.             
\end{equation}
\end{theorem}
\begin{proof}
By \ref{def:manyD} and Proposition \ref{prop:Sel-infinity}, there is $D_0$ such that $\dimp\ell\Sel{\lambda}(J_{D_0}, K) > n$. Then, by Proposition \ref{prop:Sel-infinity} applied to $J_{D_0}$, we prove the result.
\end{proof}

\begin{corollary}\label{cor:fermat}
Let $K=\ratn(\zetal)$ where $\ell$ is an odd regular prime number.
Let $F_D/K$ be the Fermat curve given by $x^\ell + y^\ell =D\, z^\ell$ where $D\in \zz$ is nonzero, and let $\Jac(F_D)/K$ be the Jacobian variety of $F_D/K$. Let $\zetal$ denote the automorphism of order $\ell$ on $F_D$ given by $x\mapsto x$, $y\mapsto y$, and $z \mapsto z\, \zetal$.  Let $\lambda$ be the endomorphism $1-[\zetal]$ on $\Jac(F_D)$ where $[\zetal]$ denotes the automorphism on $\Jac(F_D)$ induced by the automorphism on $F_D$.  Then
$ \lim\sup_D\ \dimp\ell\sellam(\Jac(F_D),K) =\infty$.
\end{corollary}
\begin{proof}
Let $C_D$ be the curve given by the homogeneous equation 
\begin{equation}\label{eq:fermat}
    D\ z^\ell = \sum_{k=0}^{(\ell-1)/2}  x^{\ell-2k}\ y^{2k}\binom{\ell}{2k}(2\ell)^{(\ell-1-2k)a}\ell\Inv
\end{equation}
where $a$ and $b$ are integers such that $a>0$ and $(\ell-1)a+\ell b =1$.
Then we have the isomorphism  $F_D \to C_D$ given by
$$
        (x:y:z) \mapsto (  x+y: (2\ell)^a( x-y ): 2^{1-b}\ell^{-b} z).
        $$

Dehomogenized with respect to $x$, the equation (\ref{eq:fermat}) is written $D\, z^\ell = f(y)$ for some monic polynomial $f(y)$ of degree $\ell-1$, defined over $\zz$.     
Note that the $K$-rational points $(\zetal^s:-\zetal:0)$ in $F_D$ for $1<s \le \ell$ are mapped to the $K$-rational points 
$$\left( 1:(2\ell)^a\cdot\frac{\zetal^s+\zetal}{\zetal^s-\zetal}:0 \right) \in C_D
                                \ \text{for } 1<s \le \ell.$$                                  
                                Thus, $f(y)$ has $(\ell-1)$ $K$-rational roots, and Corollary \ref{cor:fermat} follows from Theorem  \ref{MainTheorem3}.\end{proof}                             

\section{The General Reciprocity Laws}\label{sec:rec}
The general reciprocity law is used in this paper to show the existence of infinitely many prime numbers or ideals satisfying a set of conditions under which we are able to control the size of the Selmer groups. We recall it below.
\par
Let $K$ be a number field or a \globalfieldv\ such that $K$ contains a primitive $n$-th root of unity. 
Let $\primep$ be a place in $M_K^0$.  
For a non-archimedean place $\primep$,
we have the Hilbert norm residue symbol which is nondegenerate, and  bilinear:
\begin{equation*}
\hilbert \bullet \bullet \primep : \thegroupn{\Kp} \times \thegroupn{\Kp} \To \mu_n.
\end{equation*}
The norm residue symbol extends for an archimedean place $v$ of $K$ in an obvious way.
\par
Let $F$ be a global field.
Let $\primep$ be a prime ideal of $\OF$, and $\al \in \OF$ such that $\primep$, $\al\OF$, and $n\OF$ are pairwise coprime.
We define
\begin{equation*}
\legendren{\al}{\primep} := \zeta_n^j \text{ for some $j$ such that } \zeta_n^j \equiv \al^{ ((\Sharp \OF/\primep) -1)/n } \mod \primep.
\end{equation*}
Note that $\legendren{\al}{\primep}=1$ if and only if $x^n - \al$ has a root mod $\primep$, (provided that $\al \not\equiv 0 \mod \primep$).

If $\beta\in\OF$ is coprime to $n$ and $\al$, we extend the power residue symbol as follows: let 
$\beta\OF =\primep_1^{a_1} \cdots \primep_t^{a_t}$ be the prime ideal decomposition.
\begin{equation*}
        \legendren{\al}{\beta} := \prod_{i=1}^t \legendren{\al}{\primep_i}^{a_i}.
\end{equation*}  

\begin{definition}\label{def:power-symbol}
Let $\ell$ be a prime number. Let $K$ be a number field containing $\zetal$. 
For an \arch\ place $v\in \MK^\infty$, and for a prime ideal $\primep$ of $\OK$ which divides $\ell\OK$, we extend the symbol as follows only for convenience: Let $e:=\ord_\primep(\ell\OK)$, and let $m$ be the smallest integer which is (strictly) greater than $e\ell/(\ell-1)$.
\begin{gather*}
\legendre{\al}{v}:=
        \begin{cases}
        1 & \text{ if $x^\ell-\al=0$ is solvable over $K_v$ },\\
        -1 & \text{ if $x^\ell-\al=0$ is not solvable over $K_v$ }\\
        0 & \text{ if $\al=0$}
        \end{cases},\\
\legendre{\al}{\primep}:=
        \begin{cases}
        1 & \text{ if $\al\equiv a^\ell \mod \primep^{m}$ for some $a\not\equiv 0 \mod \primep$},\\
        -1 & \text{ if $\al\not\equiv a^\ell \mod \primep^{m}$ for all $a\not\equiv 0 \mod \primep$}\\
        0 & \text{ if $\al\equiv 0 \mod \primep$}
        \end{cases}.        
\end{gather*}
\end{definition}

\begin{theorem}\label{thm:hilbert}\ThmAuthor{The General Reciprocity Law}
Let $K$ be a number field or a \globalfieldv. Suppose that $K$ contains a primitive $n$-th root of unity where $n\ge 2$ is a positive integer coprime to $\Char K$. 
\par
Let $\al$ and $\beta$ be non-zero elements of $\OK$ coprime to each other and to $n$.
Then
$$\legendren{\al}{\beta}\cdot {\legendren{\beta}{\al}}^{-1} 
                =\prod_{\primep \in S_\infty} \hilbert\al\beta\primep$$
where $S_\infty:=\set{ v \in \MK : v \mid n\text{ or } v \in \MKinf}$.
\end{theorem}
\begin{proof} See \cite{neukirch:1999}, Theorem 8.3, p.415, or \cite{tate2:1967}, p.~352\end{proof}

\begin{lemma}\label{cor:rec1}
Let $\ell$ be a  prime number.
Let $q$ be a rational prime coprime to $\ell$, and let $K$  be the $\ell$-th cyclotomic extension of $\ratn$.
Let $\lambda:=1-\zetal$. 
Let $\primep$ be a prime ideal of $\OK$ lying over $q$.
Suppose that $\ell$ is odd. If $\legendre{\zetal}{\primep}=1$, then $q \equiv a^\ell \mod \ell^2$ for some $a\in\zz$, and  $\legendre{q}{\lambda\OK}=1$.
\par
Suppose that $\ell=2$, i.e., $K=\ratn$.  If $q$ is an odd prime such that $\legendre{-1}{q}=\legendre{2}{q}=1$, then $\legendre{q}{2\zz}=1$.
\end{lemma}
\begin{proof}
Suppose that $\ell$ is odd, and $\legendre{\zetal}{\primep}=1$.
By definition, 
$$1=\legendre{\zetal}{\primep}=\zetal^{ (q^m - 1)/\ell }
\text{ where $m$ is the residue degree of $\primep$ over $q$.}
$$
Hence, $(q^m-1)/\ell \equiv 0 \mod \ell$, i.e., $q^m\equiv 1 \mod \ell^2$. 
Note that $K/\ratn$ being Galois implies that $m \mid (\ell-1)$ and, hence, $(m,\ell)=1$.
Then there are integers $a$ and $b$ such that $am+b\ell=1$.  It follows that $q=q^{am+b\ell}\equiv (q^b)^\ell \mod \ell^2$.
By definition, 
$$ \legendre{q}{\lambda\OK}=1 \text{ if and only if }q\equiv \gamma^\ell \mod \lambda^{\ell + 1} 
                                \text{ for some } \gamma\in\OK.$$
Since $\ell\OK=\lambda^{\ell-1}\OK$, we have $q \equiv (q^b)^\ell \mod \lambda^{2(\ell-1)}\OK$.  Hence, if $\ell\ge 3$, then $2(\ell-1)\ge\ell+1$ and,  hence,$\legendre{q}{\lambda\OK}=1$.  
\par
Suppose that $\ell=2$.  Then $\zetal=-1$.  By the  supplementary quadratic reciprocity laws, 
$1=\legendre{-1}{q}=(-1)^{(q-1)/2}$ and $1=\legendre{2}{q}=(-1)^{(q^2-1)/8}$ and, hence, $q\equiv 1 \mod 4$, and $q^2\equiv 1 \mod 16$.  It follows that $q\equiv 1 \mod 8$. If $q \equiv 1 \mod 8$, then $\legendre{q}{2\zz}=1$.
\end{proof}
\begin{corollary}\label{cor:rec2}
If $a$ is a positive rational integer coprime to $\ell$ such that $\legendre{a}{\lambda\OK}=1$, then 
$$\legendre{a}{\al} = \legendre{\al}{a}$$
for all $\al\in\OK$ coprime to $a\ell$.
\end{corollary}
\begin{proof} 
Since $\legendre{a}{\lambda\OK}=1$, by definition, it implies that $a\in (K_\lambda^*)^\ell$ and, hence, the Hilbert residue symbol  $\hilbertl \al a  \lambda =1$.
If $\ell=2$, then there is an infinite place $v$ for which $K_v\cong \real$.  Since $a\in(\real^*)^2$, $\hilbertl \al a v =1$.
If $\ell\ge 3$, then the Hilbert residue symbol is trivial at all infinite places.  Therefore, 
$\prod_{\primeq \mid \ell\infty} \hilbertl \al a \primeq =1$
 and, hence,  the assertion follows from Theorem \ref{thm:hilbert}.\end{proof}
 \par
 Let $k$ be a finite field of characteristic $q$, and let $K/k$ be a \globalfieldv\ (with $\MKinf=\set{v_\infty}$).
 Let $\piinf$ be a uniformizer of the discrete valuation ring $\OK[\infty]$ of $K$ at $v_\infty$, and let $\ordinf:=\ord_{\vinf}$.
 Recall that a nonzero element $\al \in \OK$ is \textit{monic} (with respect to $\pi_\infty$) if $\piinf^m \al \equiv 1 \mod \piinf$ for some integer $m \in \zz$.

 \begin{theorem}\label{thm:rec-ff}
 \ThmAuthor{The general reciprocity law for function fields}
 Let $q$ be a prime number, and let $n$ be a positive integer not divisible by $q$.
 Let $k$ be a finite field with $q^r$ elements such that $k$ contains a primitive $n$-th root of unity.
 Let $K/k$ be a \globalfieldv. 
 \par
 If $g$ and $h$ are monic distinct elements of $\OK$ such that $g$ is coprime to $h$, then
 \begin{gather*}
 \legendren{-1}{g} = (-1)^{ \big( (q^r-1)/n\big) \cdot\ \ordinf(g)  };\\
 \rule{0pt}{1.6\baselineskip}\legendren{g}{h}\, \legendren{h}{g}\Inv = (-1)^{ \big((q^r-1)/n\big) \cdot\ \ordinf(g)\, \ordinf(h) }.
 \end{gather*}
 \end{theorem}
 \begin{proof}
 Let $\Kinf$ be the completion of $K$ at $v_\infty$. Let $\Oinfhat:=\set{\al \in \Kinf : \abs{\al}_{\vinf} \le 1}$, and 
 $\Minfhat:=\set{\al \in \Oinfhat : \abs{\al}_{\vinf} < 1}$.
 Since $\deg(\vinf)=1$, let us define $\ome : \Oinfhat^* \to k^*$ by $\al \mapsto a$ such that $\al \equiv a \mod \Minfhat$.
 By \cite{neukirch:1999}, Chapter V, Sec 3, Proposition 3.4, 
 for nonzero elements $\al$ and $\beta$ in $\OK$, 
 \begin{equation}\label{eq:omega}
 \hilbert{\al}{\beta}{\vinf} = \ome\left(     
                        (-1)^{ \ordinf(\al)\, \ordinf(\beta) } \ \frac{ \beta^{\ordinf(\al)} }{ \al^{\ordinf(\beta)} } \right)^{(q^r-1)/n}.
 \end{equation}
 Let $\piinf \in K^*$ be a uniformizer of $\Kinf$. Let $g$ and $h$ be monic distinct elements of $\OK$ coprime to each other.
 Then, by Theorem \ref{thm:hilbert},
 $$
 \legendren{-1}{g}=\hilbert{-1}{g}{\vinf} = \ome\big( (-1)^{\ordinf(g)} \big)^{ (q^r-1)/n } 
                                = (-1)^{\ordinf(g)\ (q^r-1)/n }.$$
Since $g$ and $h$ are monic, there are $a$ and $b$ in $\Oinfhat^*$ such that $a\equiv b\equiv 1 \mod \Minfhat$, $g=a\, \piinf^{ \ordinf(g) }$, and $h = b \piinf^{ \ordinf(h) }$. It follows from Hensel's lemma that $a$ and $b$ are contained in $(\Kinf^*)^n$.
Then, by Theorem \ref{thm:hilbert}, 
$$
        \legendren{g}{h}\ \legendren{h}{g}\Inv = \hilbert{\piinf}{\piinf}{\vinf}^{\ordinf(g)\, \ordinf(h)}.$$
By (\ref{eq:omega}), 
$$\hilbert{\piinf}{\piinf}{\vinf} = (-1)^{(q^r-1)/n}.$$
\end{proof}

\begin{lemma}\label{lem:rec2}
If $g$ is a monic element of $\OK$ such that $\legendren{-1}{g} =1 $, then for all monic elements $h$ of $\OK$ coprime to $g$, 
$$ \legendren{h}{g} = \legendren{g}{h}.$$
\end{lemma}
\begin{proof}
The condition: $\legendren{-1}{g}=1$ implies that $\ordinf(g)\cdot (q^r-1)/n \equiv 0 \mod 2$ and, hence,\\ $\ordinf(h)\,\ordinf(g)\cdot (q^r-1)/n \equiv 0 \mod 2$. 
By Theorem \ref{thm:rec-ff}, we proved the lemma.\end{proof}

\begin{lemma}\label{hensel-lemma} 
\par
Let $\al$ be an element of $\OK$, and let $\primep \in M_K$.  Then $\legendre{\al}{\primep}=1$ implies that $\al\in(K_\primep^*)^\ell$.
\end{lemma}
\begin{proof} 
The only non-trivial case is that $\primep \mid \ell$.
Suppose that $\primep$ is a prime ideal of $\OK$ such that $\primep \mid \ell$ and $\legendre{\al}{\primep}=1$.
Let $m$ be an integer $> e\ell/(\ell-1)$ where $e:=\ord_\primep(\ell\OK)$.
Then $\al \equiv a^\ell \mod \primep^m$ for some $a \in \OK$ implies that $\al \in \ellthpower{K_\primep}$.
\end{proof}

\begin{acknowledgements}
I wish to thank Professor Dino Lorenzini for many  useful suggestions on this work as well as much of his help on the exposition of this paper.
\end{acknowledgements}
\par

\renewcommand{\baselinestretch}{0.9}
\small
 \bibliography{thesis_pub}

\providecommand{\bysame}{\leavevmode\hbox to3em{\hrulefill}\thinspace}
\providecommand{\MR}{\relax\ifhmode\unskip\space\fi MR }
\providecommand{\MRhref}[2]{%
  \href{http://www.ams.org/mathscinet-getitem?mr=#1}{#2}
}
\providecommand{\href}[2]{#2}
\begin{thebibliography}{{}Neu99}

\bibitem[Ata01]{atake:2001}
D.~Atake, \emph{{On elliptic curves with large Tate-Shafarevich Groups}}, J.\
  Number Theory \textbf{87} (2001), 282--300.

\bibitem[Cha]{chang:2004}
S.~Chang, \emph{{ Note on the rank of quadratic twists of Mordell equations,
  \rm to appear in J. Number Theory}}.

\bibitem[CJB97]{chm:1997}
L.~Caporaso, {J. Harris}, and {B. Mazur}, \emph{{Uniformity of rational
  points}}, J.~Amer.~Math.~Soc. \textbf{10} (1997), 1--35.

\bibitem[HB94]{heath-brown:1994}
D.R. Heath-Brown, \emph{{ The size of Selmer groups for the congruent number
  problem II}}, Invent.\ Math. \textbf{118} (1994), 331--370.

\bibitem[IP00]{iwaniec-sarnak:2000}
H.~{}Iwaniec and {P. Sarnak}, \emph{{The non-vanishing of central values of
  automorphic $L$-functions and Landau-Siegel zeros}}, Israel J.\ Math
  \textbf{120} (2000), 155--177.

\bibitem[{}Kol88]{kolyvagin:1988}
V.A. {}Kolyvagin, \emph{{Finiteness of $E(\ratn)$ and $\sha(E,\ratn)$ for a
  subclass of Weil curves}}, Izv.~Akad.~Nauk SSSR Serl.~Mat. \textbf{52}
  (1988), 1154--1180.

\bibitem[Lem]{lemm:1998}
F.~Lemmermeyer, \emph{{ On Tate-Shafarevich groups of some elliptic curves}},
  in \lq\lq Proc. Conf., Graz\rq\rq, 1998.

\bibitem[{}Lie94]{lieman:1994}
D.~{}Lieman, \emph{{Non-vanishing of $L$-series associated to cubic twists of
  elliptic curves}}, Ann.\ of Math. \textbf{140} (1994), 181--108.

\bibitem[{}Lor96]{dino:1996}
D.~{}Lorenzini, \emph{{Invitation to Arithmetic Geometry}}, AMS, 1996.

\bibitem[LT02]{dino:2002}
D.~{}Lorenzini and {T. Tucker}, \emph{{Thue equations and the method of
  Chabauty-Coleman}}, Invent.\ Math. \textbf{148} (2002), 47--77.

\bibitem[Maz86]{mazur:1986}
B.~Mazur, \emph{{Arithmetic on curves}}, Bull.\ Amer.\ Math.\ Soc. \textbf{14}
  (1986), 207--260.

\bibitem[{}Neu99]{neukirch:1999}
J.~{}Neukirch, \emph{{Algebraic Number Theory}}, Springer-Verlag Berlin
  Heidelberg, 1999.

\bibitem[OC98]{ono:1998}
K.~{}Ono and {C. Skinner}, \emph{{ Non-vanishing of quadratic twists of modular
  $L$-functions}}, Invent.\ Math. \textbf{134} (1998), 651--660.

\bibitem[PE97]{poonen:1997}
B.~Poonen and {E. Schaefer}, \emph{{Explicit descent for Jacobians of cyclic
  covers of the projective line}}, J.\ reine angew.\ Math. \textbf{488} (1997),
  141--188.

\bibitem[{}Sch90]{schoen:1990}
C.~{}Schoen, \emph{{Bounds for rational points on twists of constant
  hyperelliptic curves}}, J.\ reine angew.\ Math. \textbf{411} (1990),
  196--204.

\bibitem[Sch98]{schaefer:1998}
E~Schaefer, \emph{{Computing a Selmer group of a Jacobian using functions on
  the curve}}, Math.~Ann. \textbf{310} (1998), 447--471.

\bibitem[{}Ser76]{serre:1976}
J.P. {}Serre, \emph{{Divisibilit\'e de certaines fonctions arithm\'etiques, \rm
  Enseign. Math. \bf 22\rm\ (1976), 227--260}}, Enseign.\ Math. \textbf{22}
  (1976), 227--260.

\bibitem[{}Sil93]{silverman:1993}
J.~{}Silverman, \emph{{A uniform bound for rational points on twists of a given
  curve}}, J.~London Math.~Soc (2) \textbf{47} (1993), 385--394.

\bibitem[{}Sto]{stoll:2004}
M.~{}Stoll, \emph{{Independence of rational points on twists of a given curve,
  \rm submitted}}.

\bibitem[{}Sto98]{stoll:1998}
\bysame, \emph{{On the arithmetic of the curves $y^2=x^\ell+A$ and their
  Jacobians}}, J. reine angew. Math. \textbf{501} (1998), 171--189.

\bibitem[Tat]{tate2:1967}
J.~Tate, \emph{{Fourier analysis in number fields and Hecke's
  Zeta-functions,\rm\ In: J.W.S.~Cassels and A.~Fr\"ohlich (eds): Algebraic
  Number Theory, Academic Press 1967 }}.

\bibitem[{}Vat98]{vatsal:1998}
V.~{}Vatsal, \emph{{Rank-one twists of a certain elliptic curve}}, Math.~Ann.
  \textbf{311} (1998), 791--794.

\bibitem[Won99]{wong:1999}
S.~Wong, \emph{{Elliptic curves and class number divisibility}}, Internat.\
  Math.\ Res.\ Notices \textbf{12} (1999), 661--672.

\bibitem[{}Yu03]{yu:2003}
G.~{}Yu, \emph{{Rank $0$ quadratic twists of a family of elliptic curves}},
  Compositio Math. \textbf{135} (2003), 331--356.

\end{thebibliography}
 \bibliographystyle{amsalpha}

\end{document}